\documentclass[11pt]{article}

\usepackage{mathrsfs,amssymb,amsmath,amsthm}
\usepackage{authblk}

\usepackage{graphicx}
\usepackage[dvipsnames]{xcolor}
\usepackage{anyfontsize}
\usepackage{tikz}
\usetikzlibrary{shapes,positioning}
\usepackage[top=50pt, bottom=80pt,left=70pt, right=70pt]{geometry}
\usepackage{soul}

\usepackage[pagebackref,colorlinks,citecolor=Plum,urlcolor=Periwinkle,linkcolor=DarkOrchid]{hyperref}
\usepackage[square,numbers]{natbib}
\usepackage{verbatim} 
\RequirePackage[labelsep=space,tableposition=top]{caption}

\usepackage{subcaption}
\usepackage[left,pagewise]{lineno} 

\newcommand{\bb}[1]{\mathbb{#1}}

\newcommand{\pv}[2]{#1 a^\perp \otimes a^\perp + #2 n^\perp \otimes n^\perp}
\newcommand{\pva}[1]{#1 a^\perp \otimes a^\perp}
\newcommand{\pvn}[1]{#1 n^\perp \otimes n^\perp}

\newcommand{\re}{\ensuremath{\mathbb{R}}}
\newcommand{\ten}[2]{#1_{#2} \otimes #1_{#2}}
\newcommand{\sten}[2]{#1 \odot #2}
\newcommand{\m}{\re^{2\times 2}}
\newcommand{\msym}{\re^{2\times 2}_{sym}}
\newcommand{\mskew}{\re^{2\times 2}_{skew}}
\newcommand{\pro}[2]{\left \langle #1 \, , #2 \right \rangle}

\DeclareMathOperator{\rpartial}{\partial_{ri}}
\DeclareMathOperator{\Aff}{Aff}
\DeclareMathOperator{\rank}{Rank}

\DeclareMathOperator{\Ker}{Ker}
\DeclareMathOperator{\kers}{Ker_s}
\DeclareMathOperator{\sgn}{sign}
\DeclareMathOperator{\Id}{Id}

\newtheorem{definition}{Definition}[section]
\newtheorem{theorem}{Theorem}[section]
\newtheorem{proposition}{Proposition}[section]
\newtheorem{corollary}{Corollary}[theorem]
\newtheorem{remark}{Remark}[section]
\newtheorem{lemma}{Lemma}[section]

\newtheorem*{theorem*}{Theorem}
\newtheorem*{proposition*}{Proposition}
\newtheorem*{remark*}{Remark}
\newtheorem*{lemma*}{Lemma}

\everymath{\displaystyle}

\definecolor{silver}{rgb}{0.75, 0.75, 0.75}\definecolor{trolleygrey}{rgb}{0.5, 0.5, 0.5}
\definecolor{pansypurple}{rgb}{0.0, 0.13, 0.29}
\definecolor{mirojo}{RGB}{200,4,90}
\definecolor{miverde}{RGB}{4,220,100}
\definecolor{mypurple}{RGB}{150,10,120}
\definecolor{morado}{RGB}{125,0,177}
\definecolor{rojo}{RGB}{255,0,0}
\definecolor{azul}{RGB}{0,0,255}
\definecolor{verde}{RGB}{0,135,0}

\usepackage[disable]{todonotes} 
\newcommand{\notelm}[1]{\todo[linecolor=green,backgroundcolor=green!40]{#1}}

\makeatletter
\if@todonotes@disabled

\else

\fi
\makeatother

\makeatletter
\if@todonotes@disabled
\newcommand{\todoin}[1]{#1}
\else
\newcommand{\todoin}[1]{\todo[inline, caption={2do}, backgroundcolor=red!20, bordercolor=red!20, size=\normalsize]{
 \begin{minipage}{\textwidth-8pt}#1\end{minipage}}}
\fi
\makeatother

\makeatletter
\if@todonotes@disabled

\else

\fi
\makeatother

\definecolor{morado}{RGB}{125,0,177}
\newcommand{\kca}{\color{black}} 

\usepackage{float}
\restylefloat{figure}

\usepackage{booktabs} 
\usepackage{multirow}

\usepackage{cleveref} 

\usepackage{enumitem}
\begin{document}

\date{December 2020}

\author[1]{Capella, Antonio}
\author[2]{Morales, Lauro}
\affil[1]{\small Instituto de Matem\'aticas-UNAM \\ lmm@ciencias.unam.mx
}
\affil[2]{\small Instituto de Matem\'aticas-UNAM\\
capella@im.unam.mx}

\title{On the quasiconvex hull for a three-well problem in two dimensional linear elasticity}
\maketitle
\begin{abstract}
\hskip -.2in
\noindent
We provide quantitative inner and outer bounds for the symmetric quasiconvex hull $Q^e(\mathcal{U})$ on linear strains generated by three-well sets $\mathcal{U}$ in $\msym$. In our study, we consider all possible compatible configurations for three wells and prove that if there exist two matrices in  $\mathcal{U}$ that are rank-one compatible then $Q^e(\mathcal{U})$ coincides with its symmetric 
lamination convex hull $L^e(\mathcal{U})$. We complete this result by providing an explicit characterization of $L^e(\mathcal{U})$ in terms of the wells in  $\mathcal{U}$. 
Finally, we discuss the optimality of our outer bound and its relationship with quadratic polyconvex functions. 
\end{abstract}

\section{Introduction}\label{intro}

In this work, we are concerned with variational problems of the form 
\begin{equation}\label{freeEC}
\min_{y=y_0 \mbox{ on } \partial \Omega} I(y),\quad \mbox{where} \quad I(y)=\int_\Omega f(\nabla y(x))dx,
\end{equation}
the function $y:\Omega\subset\re^n \to \re^n$, and $f:\re^{n\times n}\to \re$ is continuous, bounded from below, and satisfies mild growth conditions (see, e.g. section 9 in \cite{BDa}). 
For integral functionals like \eqref{freeEC}, other properties of $f$ besides continuity and coercivity are needed to assure the sequential weak lower semicontinuity (swlsc) of $I$, {\em i.e.}, the weak convergence of minimizing sequences to minimizers, see  for example theorems 1.13 and 1.15 in \cite{BDa} pp. 15 and 17, and section 3 in \cite{BK}. The appropriate additional condition is some kind of convexity of $f$.
In the 1950s, Morrey, see theorems 2.1 and 2.2 in \cite{CM} pp. 5 and 8, respectively,  introduced quasiconvexity as the necessary and sufficient condition for swlsc (see also theorem 3.4 and remark 3.3 in \cite{BK}).  
A locally Borel-measurable function $f$ is {\it quasiconvex} if for every $M\in \re^{n\times n}$ and each smooth function $\phi:\Omega \to \re$ 
compactly supported in $\Omega$, 
\begin{equation*}
f(M) \leq \frac{1}{|\Omega|}\int_\Omega f(M+ \nabla \phi)\,dx.
 \end{equation*}

In materials science, energy models like \eqref{freeEC} arise in the study of martensitic phase transformations. 
Such models satisfy frame indifference and multi-well structure assumptions. 
By multi-well structure, we mean that $f\geq 0$ everywhere but $f=0$ only on a prescribed set $K$. The connected components of $K$ are the energy wells that represent the material's different phases. 
In most cases, these energy densities are not quasiconvex under reasonable conditions, see  \cite{BJ2F} p. 401, and  \cite{K} p. 1.  Thus, minimizing sequences may develop oscillations on their gradient, and only weak convergence is to be expected. In these models, oscillations correspond to the microstructures observed in real materials, see \cite{BJ} p. 14. 
There are three types of sets $K$ relevant to models of martensitic materials, 
 depending on $f$'s symmetry. First, $f$ has no symmetry and $K\subset \re^{n\times n}$ is a set of {\em matrices}. Second, $f$ is left-invariant by orthogonal matrices,  i.e., $f(OM)=f(M)$ for every $M\in \re^{n\times n}$ and  $O\in SO(n)$. In this case, $K =\{SO(n) M_1, ..., SO(n)M_k\}$, where $SO(n)M =\{Q M \ | \ Q\in SO(n)\}$ is the elastic well associated with the symmetric positive-defined matrix $M$.  Third, $f$ is invariant under the addition of skew-symmetric matrices. In this case, $K=\{M_1\oplus \re^{n\times n}_{skew},...,M_k\oplus \re^{n\times n}_{skew}\}$ is a set of {\em linearized elastic wells} or {\em linear strains}, where  $M\oplus \re^{n\times n}_{skew}:=\{M+W\ |\ W\in \re^{n\times n}_{skew} \}$ is the linearized elastic well of $M$. We also say that $M$ represents the well $M\oplus \re^{n\times n}_{skew}$.
 Elements in $K$ can be used to construct minimizing sequences in terms of sequential lamination
 provided they are compatible, see for example \cite{K} pp. 219 -- 220 and theorem 3 in \cite{BJ2F} pp. 25 -- 26.
As in  \cite{BFJK} p. 846, we say that two {\em matrices} $M_1,\,M_2\in \re^{n\times n}$, are compatible if they are {\em rank-one connected},  i.e., $\rank(M_1-M_2) \leq 1$. 
In the opposite case they are called {\em incompatible}.  
Two elastic wells $SO(n)M_1$ and  $SO(n)M_2$, are {\em compatible}
if there exists $Q_1,Q_2\in SO(n)$ such that $Q_1M_1$ and $Q_2M_2$ are rank-one connected. They are called incompatible otherwise. 
The corresponding notion of compatibility  in linear elasticity is that: two linear strains  $M_1\oplus \re^{n\times n}_{sym}$ and $M_2\oplus \re^{n\times n}_{sym}$ are compatible if there exists a skew symmetric matrix $W$  such that   $\rank(M_1-M_2+W) \le 1$. In the particular case where $M_2=0$, we simply say that $M_1$ is compatible or incompatible if $\det M_1\leq 0$ of  $\det M_1>0$, respectively.  As before, they are called incompatible otherwise.  
Moreover, we say that two matrices $M_1,M_2\in\msym$ which represent linear strains are {\em rank-one compatible} if $\det (M_1-M_2)=0$. 
 Moreover,  if $M_1$ and $M_2$ represent two  compatible linear strains,  then $M_1-M_2+W=a\otimes n$  
 for some  $W\in \re^{2\times 2}_{skew}$,  $a\in \re^2$, and $n\in \mathbb{S}^1$. 

If two matrices, elastic wells or linear strains, are compatible,  layering structures can be constructed and  nontrivial oscillating minimizing sequences emerge in suitable spaces, see for example proposition 2 in \cite{BJ} p. 23, and section 3 in \cite{BFJK} pp. 855 to 868. 
Young measures are an effective tool to capture the nonlinear functional's asymptotic behavior along an oscillating sequence.  While the weak limit carries information on average values in an oscillating sequence, the Young measure contains information on ``where" these oscillations occur.  
If the energy achieves its minimum, the minimizing sequence does not oscillate,  and the corresponding Young measure is trivial; 
namely, it is a Dirac mass for almost every point in the domain. 
For one-well and two-incompatible-wells problems, Kinderleherer \cite{Kin} (pp. 15 and 16) showed that the Young measure limit of gradients (case of  elastic wells above) is trivial and furthermore constant in the domain. The case for two linear strains was studied by Kohn in \cite{K}.
Bhattacharya et al. \cite{BFJK},  studied the problem of which measures arise as Young measures limits of  gradients, linear strains, and matrices. In the case of $k$ pairwise incompatible elastic wells or linear strains in two dimensions (see also Lemma 3 in 
\cite{Sv3}, p. 408), they show that the corresponding Young measures' limits are trivial.  
The situation changes as soon as we consider four matrices in two space dimensions, see remark 6 in \cite{Tar} pp. 194 and 195, three elastic wells, or three linear strains in three dimensions.  The authors in \cite{BFJK} showed that in the later three cases there exist nontrivial Young measures in the limit under the assumption of pairwise incompatibility. They called this phenomenon 
``mutual compatibility", see  \cite{BFJK} p. 849.

Depending on the space dimension and the number of compatible --or incompatible-- pairs among the matrices, wells, or linear strains in $K$, we may have the existence of nontrivial Young measures in the limit. 
A closely related problem is to identify the set of values taken by the (constant) weak limits of minimizing sequences.  
This set is known as the quasiconvex hull $Q(K)$ of $K$. 
If the corresponding Young measures are trivial, then $Q(K)=K$. As stated before, that is the case for any number of pairwise incompatible linear strains in two dimensions. If there exists a nontrivial Young measure in the limit, as in the above examples,  we have that $K$ is strictly contained in $Q(K)$. In the fully pairwise compatible case, Bhattacharya~\cite{B} p. 231 --see also Lemma 4 in \cite{A2016} p.46 -- proved that the quasiconvex and convex hulls are equal for any number linear strains in $K$.  

Quasiconvex hulls are generically very difficult (if not impossible) to compute for a particular choice of $\mathcal{U}$.   
Explicit examples of quasiconvex functions \cite{BDa,Sv4}  and quasiconvex hulls \cite{BD, BFJK, Bt, Z1998} are scarce and most  of the nontrivial examples are not known explicitly. 
One approach to this problem is to compute the quasiconvex envelope of the energy density and to study its zero energy level set. This is a difficult problem, and solutions are only known in few cases, for example  see theorems 3.4 and 3.5 in \cite{K} pp 204 and 205. 
There is also extensive literature on the exploration of upper and lower bounds for the quasiconvex envelope of a function, see for example \cite{Dol,KL,Anja2015}.

An alternative approach is to consider quasiconvex hulls' inner and outer bounds, such as the polyconvex, rank-one, and lamination convex hulls.   
Following \cite{BDa} pp. 156 and 157, we introduce some standard notions from semi-convex analysis. We say that a function $f:\re^{n\times n} \rightarrow \re$ is {\em polyconvex} 
if there exists a convex function $G:\re^{m_n}\rightarrow \re$, such that 
$f(M) = G\circ T(M)$,
where $T:\re^{n\times n}\rightarrow \re^{m_n}$  is given by  
\begin{equation}\label{Tadj}
T(M) = (M, adj_2(M),\ adj_3(M),\ \cdots, \det M ),
\end{equation}
$adj_k(M)$ stands for the $k\times k$ matrix of sub-determinants of $M$,  $m_n=C(2n,n)-1,$ and 
$C(2n,n)$ is the binomial coefficient between $2n$ and $n$. 
The function $f:\re^{n\times n} \rightarrow \re$ is {\em rank-one convex} if for every $M_1,M_2\in \re^{n\times n}$ such that $\rank \, (M_1-M_2)\leq 1$, 
and every $\lambda\in [0,1]$,
\[
f(\lambda M_1 +(1-\lambda)M_2) \leq \lambda f(M_1) +(1-\lambda)f(M_2).
\]
It can be proved that if $f$ is convex, it is polyconvex, if $f$ is polyconvex, it is quasiconvex, and if $f$ is quasiconvex, it is rank-one convex,
but the reverse implications are false in general, see sections 5.3.2 to 5.3.9 in  \cite{BDa} and \cite{Sv2}. 

For any compact set $K\subset \re^{n\times n}$, its semi-convex hulls are defined by means of cosets, more precisely
\[
S(K) = \left\{A\in \re^{n\times n}\, | \, f(A) \leq \sup_{B\in K} f(B), \ \mbox{for every semi-convex}  \ f  \right\}.
\] 
It is clear that $R(K)\subset Q(K)\subset P(K)\subset C(K)$  where these sets  correspond to the rank-one, 
quasiconvex, polyconvex and convex hulls, respectively. These inclusions may not be proper, see for example Theorems 7.7 and 7.28 in \cite{BDa}.


In the geometrically linear regime, 
 the relevant quantities \cite{K} are given by the {\it displacement } $u(x)=y(x)-x$ and the {\it linear strain}
\begin{equation*}
e(\nabla u):=\frac{\nabla u + (\nabla u)^T}{2}.
\end{equation*}
In this framework, the energy density becomes invariant under  the addition of skew-symmetric matrices in its argument, {\em i.e.}, $f(M)=f(M+S)$ for every  $ S\in \re^{n\times n}_{skew}$, 
and it vanishes on the set of linearized elastic wells  represented by $\mathcal{U} = \{U_1,\, U_2,\, \cdots U_k\}\in \re^{n\times n}_{sym}$.
We use the notation $K=\mathcal{U}\, \oplus\, \re^{n\times n}_{skew}$  to indicate 
that the symmetric part of every element in $K$ belongs to $\mathcal{U}$. 
Also, similar ideas about semi-convexity are available.
Following Boussaid {\em et. al} \cite{Anja19} (section 2, pp 423 to 427), 
we say that a function $f:\re^{n\times n}_{sym} \rightarrow \re$ is symmetric semi-convex if $f(e(\cdot)):\re^{n\times n} \rightarrow \re$ is semi-convex.
\kca
It is straightforward to see that $f:\re^{n\times n}_{sym} \rightarrow \re$ is  symmetric quasiconvex if and only if 
for every $U\in \re^{n\times n}_{sym}$ 
\[
f(U)\leq \inf\left\{ \frac{1}{|\Omega|}\int_\Omega f(U +e(D\phi))dx \, \middle | \, \phi\in C^\infty_0(\Omega, \re^{n} ) \right\},
\] 
and it is  symmetric rank-one convex if and only if for every $U_1$ and $U_2$ compatible  symmetric matrices and $\lambda\in [0,1]$,
\[
f(\lambda U_1 +(1-\lambda)U_2) \leq \lambda f(U_1) +(1-\lambda)f(U_2). 
\]
Characterization of symmetric polyconvexity for any dimension is rather challenging. Boussaid {\em et. al} \cite{Anja19} gave an explicit characterization of symmetric polyconvex functions in the cases of two and three dimensions. 
In the two-dimensional case, they proved that $f:\msym\rightarrow \re$ is symmetric polyconvex if and only if there exists 
$g:\msym \times \re \rightarrow \re$ convex,  non-increasing on its second argument, and such that 
$f(\cdot) = g(\cdot, \det(\cdot))$.
\medskip 

Also, the symmetric semi-convex hull, see \cite{Z2002} p. 567, for a compact set $\mathcal{U}\in \re^{n\times n}_{sym}$ is 
defined in terms of the cosets, 
\begin{equation}\label{eq:cosets}
S^e(\mathcal{U}) = \left\{A\in \re^{n\times n}\, | \, f(A) \leq \sup_{B\in \mathcal{U}} f(B), \ f \ \mbox{symmetric semi-convex}\right\},
\end{equation}
where $S$ must be replaced by $R,Q$ and $P$ for the symmetric rank-one, quasiconvex and polyconvex hull of the set $\mathcal{U}$.
As in the nonlinear case, to determine the symmetric quasiconvex hull of a compact set $\mathcal{U}$ is a challenging task. 
An inner approximation for the symmetric quasiconvex hull is given by $L^e(\mathcal{U})$, the symmetric lamination convex hull of $\mathcal{U}$, see \cite{Z2002}. This set is defined  as the intersection of the symmetric lamination hulls of all ranks, namely
\[
L^e(\mathcal{U}) = \bigcup_{i=0}^\infty L^{e,i}(\mathcal{U}),
\]
where $L^{e,0}(\mathcal{U}) = \mathcal{U},$ and 
\[
L^{e,i+1}(\mathcal{U}) = \{\lambda A + (1-\lambda) B \in\re^{n\times n}_{sym} \, | \, \lambda \in [0,1]\mbox{ and } A,B\in L^{e,i}(\mathcal{U}) \mbox{ are compatible}\}. 
\]
It is also known, see for example \cite{Z2002} p. 562, that for any compact set $\mathcal{U}\subset \re^{n\times n}_{sym}$,
\begin{equation}\label{hullscontentions}
\mathcal{U} \subset L^{e}(\mathcal{U}) \subset R^e(\mathcal{U}) \subset Q^e(\mathcal{U}) \subset P^e(\mathcal{U}) \subset C(\mathcal{U}). 
\end{equation}

In this paper, we aim to understand the effect of partial pairwise  compatibility among the elements of $\mathcal{U}$ on the algebraic restrictions on the possible values of effective linear strains.
We consider sets $\mathcal{U}$ of three linearized elastic wells in two dimensions, where the wells are not fully compatible nor incompatible.  
Our goal is to characterize --or bound-- the symmetric quasiconvex hull $Q^e(\mathcal{U})$ of $\mathcal{U}$. 

Smyshlyaev and Willis \cite{SW} studied the relaxation of the three-well energy in linear elasticity where the elastic wells share the same elastic modulus. In particular, they showed that the minimizing measure may be chosen as the sum of no more than three Dirac masses. 
In \cite{Anja2015}, Schl\"omerkempe {\em et. al.} considered an analogous three-well setting in the three dimensional linear elasticity for a particular set of wells.   
Here, we do not consider any restriction on the wells to get general inner and outer bounds on $Q^e(\mathcal{U})$. For the interior bound, we characterize $L^e(\mathcal{U})$ explicitly, while for the outer bound, we estimate the zero level set of a  suitable quadratic function's polyconvex envelope.
 In some cases we are able to show that our bounds are optimal.   
\medskip

In the rest of the paper we will follow standard notation. 
For every matrix $M\in \re^{2\times 2}$,  $e(M)$ and $w(M)$ denote the symmetric and 
skew-symmetric part of $M$ respectively. Also for every $a,b\in \bb{R}^2$, the tensor product 
$a\otimes b \in \re^{2\times 2}$  is defined as $\left(a\otimes b\right)_{ij} = a_i b_j$ 
and its symmetric part is denoted by $\sten{a}{b}$. The $2\times 2$ identity matrix is denoted by $\Id$, and the bilinear form $\pro{\cdot}{\cdot}:\re^{2\times 2} \times \re^{2\times 2} \rightarrow \bb{R}$ stands for the Frobenius inner product {\em i.e }  $\pro{A}{B}\mapsto \emph{Tr} (A^T \,B)$, also $ \Vert\cdot\Vert$ denotes the Frobenius norm. 
Since the space of skew-symmetric matrices has dimension one, we will use the following representation, 
\begin{equation}\label{Rmatrix}
w(M) = w_M R, \quad \mbox{where} \quad
R = \begin{pmatrix}
0 & -1 \\ 1 & 0 \end{pmatrix} \quad \mbox{and} \quad w_M=\frac{1}{2} \pro{M}{R}.
\end{equation}
Along this work and for the sake of brevity, we shall say {\em wells } when we actually mean {\em linearized elastic wells}. Moreover, we refer to the elements of $\mathcal{U}$ as {\em wells} or as the representative of wells equally.
 
\section{Main results}\label{sec:1}

In this paper, we are interested in the  quasiconvex hull $Q^e(\mathcal{U})$ of the set  of three elastic wells  $\mathcal{U}=\{U_1,U_2,U_3\}\subset \msym$. 

\begin{theorem}\label{T3}
Let $\mathcal{U}\subset \msym$ represents a three-well set such that $\Aff( \mathcal{U})$ has codimension one. If there exist at least two different matrices in $\mathcal{U}$ that represent two rank-one compatible wells, then 
\[Q^e(\mathcal{U}) = L^e(\mathcal{U}).\] 
\end{theorem}

This result can be described in terms of the compatibility  between the wells. Independently of the existence of a rank-one compatibility,  
Bhattacharya~\cite{B} p. 231 --also see \cite{A2016} p. 46-- proved that  $Q^e(\mathcal{U}) =L^e(\mathcal{U})=C(\mathcal{U})$  if the wells are all pairwise compatible.  
On the other hand, 
Bhattacharya {\em et. al.} \cite{BFJK} (theorem 2.3 pp. 854)  showed that if
the wells are pairwise incompatible, then  $Q^e(\mathcal{U})=L^e(\mathcal{U})=\mathcal{U}$  and no 
microstructure can be formed.  The novelty of \cref{T3} is to consider the intermediate cases 
when only one or two  pairs of matrices in $\mathcal{U}$ are compatible, and one of these compatibility relation is 
a rank-one compatibility. 
\medskip

In the following result, we give an explicit characterization of the symmetric lamination convex hull when at least two wells in ${\mathcal U}$ are  incompatible.  Notice that --see \cref{qconnectness} below-- 
two wells $U_1$ and $U_2$ are compatible if and only if $\det(U_1-U_2)\leq 0$, and they are incompatible if 
$\det(U_1-U_2)>0$. Conditions $\det(U_1-U_2)=0$ and $U_1\neq U_2$ are equivalent to rank-one compatibility between $U_1$ and $U_2$.
 
 \begin{proposition}\label{P:lamconv}
Let $ \mathcal{U} = \{U_1,U_2,U_3\}$ be such that $\Aff(\mathcal{U})$ has codimension one.
\begin{enumerate}[label=(\alph*)]
\setlength\itemsep{.2em}
\item\label{lhull1}  If $\det(U_1-U_2) > 0$,  $\det(U_1-U_3)> 0$ and $\det(U_2-U_3) \leq 0$, then
\[
L^e(\mathcal{U})=  \{U_1\}  \cup C(\{U_2,U_3\}). \]
\item\label{lhull2}  If  $\det(U_1-U_2)\leq 0$, $\det(U_1-U_3)\leq 0$ and $\det(U_2-U_3)> 0$, then
\[
L^e(\mathcal{U})=  C(\{U_0,U_1,U_2\}) \cup C(\{U_0,U_1,U_3\}),
\]
where  $U_0\in C(\mathcal{U})$ is uniquely characterized by  $\det(U_0-U_3)=\det(U_0-U_2)=0$.
\end{enumerate}
  
\end{proposition}

Our final result is about the symmetric quasiconvex hull when 
there exist one incompatible and one compatible pairs of wells in
$\mathcal{U}$, and there are no rank-one compatibilities among its elements.  
With this result, we cover all possibilities  for three wells in two dimensions. \cref{T2} does not give a  complete characterization of  $Q^e(\mathcal{U})$, but it provides explicit quantitative  outer bounds. 
\begin{theorem}\label{T2}
Let $ \mathcal{U} = \{U_1,U_2,U_3\}$ represent a three-well set and assume that $\Aff(\mathcal{U})$ has codimension one. Also, assume that  (a) $\det(U_1-U_2) > 0$, $\det(U_1-U_3)> 0$, and $\det(U_2-U_3) < 0$,  or (b) $\det(U_1-U_2) <  0$, $\det(U_1-U_3)<0$, and $\det(U_2-U_3) > 0$, then  
\[Q^e(\mathcal{U}) \subset \{ U \in C(\mathcal{U}) \,| \, 0\leq \hbar(U-U_0)\} \varsubsetneq C(\mathcal{U}),\]
where $U_0\in C(\mathcal{U})$ is uniquely characterized by $\det(U_0-U_2)=\det(U_0-U_3)=0$,   $\hbar:\msym \rightarrow \re$ is given by 
\[
 \hbar(V) = \pro{C}{V-V_2}\det V_1-\pro{C}{V_1-V_2}\det V,
 \]
and 
 \[
C=
\frac{(V_2-V_1)\pro{V_3-V_1}{V_3-V_2}+(V_3-V_1)\pro{V_2-V_1}{V_2-V_3}}
{\sqrt{\|V_2-V_1\|^2\|V_3-V_1\|^2-\pro{V_2-V_1}{V_3-V_1}}}.\] 
\end{theorem}

The outer bound of \cref{T2} is not necessarily optimal. Nonetheless, we show that the outer bound in \cref{T2} coincides with the quadratic polyconvex hull of $\mathcal{U}$ at the end of this paper. 
A similar bound is given by Tang \& Zhang \cite{TZ2008} (theorem 1, p. 1266) for the case of a three-dimensional three-well problem at an exposed incompatible edge of $C(\mathcal{U})$.
In \cite{TZ2008}, Tang \& Zhang also  proved that one can chip a wedge-like slice off $C(\mathcal{U})$ at the exposed incompatible edge without touching $Q^e(\mathcal{U})$. This later estimate is independent on the diameter of the set $\mathcal{U}$, but it is neither explicit nor optimal. 
 In \cite{Anja2015}, the authors estimate by numerical approximation the quasiconvex hull for special combinations of three wells in three dimensional linear elasticity that come
 from the cubic-to-monoclinic martensitic phase transformations model. The estimated configurations are analogous to our two-dimensional results in \cref{T2} and \cref{T3}. For comparison, see cases 3 and 4 and figures 5 to 7 in \cite{Anja2015} and \cref{fig:1} below.     
 
 An important ingredient in our proofs is the geometric characterization of the subset of compatible linear strains with the null matrix.   By the isomorphism between $\re^3$ and $\msym$ (see \cref{isomorphism} below), this set can be identified with the exterior of a solid cone $\mathcal{C}_0$, whose vertex is the null matrix, and its axis is $\{t\,\Id \, |\, t\in \re\}$.  Moreover, the set  $\partial \mathcal{C}_0$ consists of rank-one symmetric matrices. The set of compatible linear strains with respect to a given matrix $U\in\msym$ can be obtained by translation of  $\mathcal{C}_0$, namely $\mathcal{C}_{U}=U+\mathcal{C}_0$, (see Lemma~\ref{conenonconnection} and Figure~\ref{fig:coneR3}).
 \medskip 
 
Theorems \ref{T3} and \ref{T2} can be geometrically interpreted by identifying $\msym$ with $\bb{R}^3$ via
the linear isomorphism
\begin{equation}\label{isomorphism}
M=\left(\begin{array}{cc} x& z \\ z & y\end{array}\right) \mapsto 
\tilde{M}=\left(\begin{array}{c} x\\y\\ z\sqrt{2}\end{array} \right).
\end{equation}

If $\mathcal{U}=\{U_1,U_2,U_3\}\subset \msym$ satisfies that $\Aff(\mathcal{U})$ has codimension one, then $\Aff(\{\tilde{U}_1,\tilde{U}_2,\tilde{U}_3\})$ defines a plane in $\re^3$. Let $U_0\in C(\mathcal{U})$ as in \cref{origin}, then 
the intersection between $\mathcal{C}_{U_0}$ and $\Aff(\mathcal{U})$ determines the set of incompatible wells with $U_0$ in $\Aff(\mathcal{U})$. Moreover, $\Aff(\mathcal{U})\cap \partial \mathcal{C}_{U_0}$ consists of matrices $U$ such that $U$ and $U_0$ represent rank-one compatible wells, and due to the isomorphism \eqref{isomorphism},
this set is identified with a pair of lines in $\re^3$. Thus we shall say that $\Aff(\mathcal{U})\cap \partial \mathcal{C}_{U_0}$ is a pair of {\em rank-one lines}, see \cref{fig:1}.
Additionally, if we assume that $\mathcal{U}$ meets the conditions in Proposition \ref{P:lamconv}, then  the lamination convex hull is explicitly known and can be represented in the plane  $\Aff(\{\tilde{U}_1,\tilde{U}_2,\tilde{U}_3\})$.  
Figures \ref{fig:1}{\color{DarkOrchid}.(a)} and \ref{fig:1}{\color{DarkOrchid}.(c)} display two configurations with two incompatibility relations between $\mathcal{U}$'s wells.  Due to Proposition \ref{P:lamconv}, the 
lamination convex hull equals the union between $\{U_1\}$ and the segment joining $U_2$ and $U_3$. In the configuration presented in figure \ref{fig:1}{\color{DarkOrchid}.(a)}, the outer bound given by \cref{T2} is the union of $\{U_1\}$ with the region bounded by the curve $\Gamma=\{U\in C(\mathcal{U}) \, | \,  \hbar(U-U_0)=0\}$ and the segment joining $U_2$ and $U_3$. Additionally, the configuration shown in  \cref{fig:1}{\color{DarkOrchid}.(c)} has a rank-one compatibility, and the set $Q^e(\mathcal{U})$ equals $L^e(\mathcal{U})$ by Theorem \ref{T3}. 
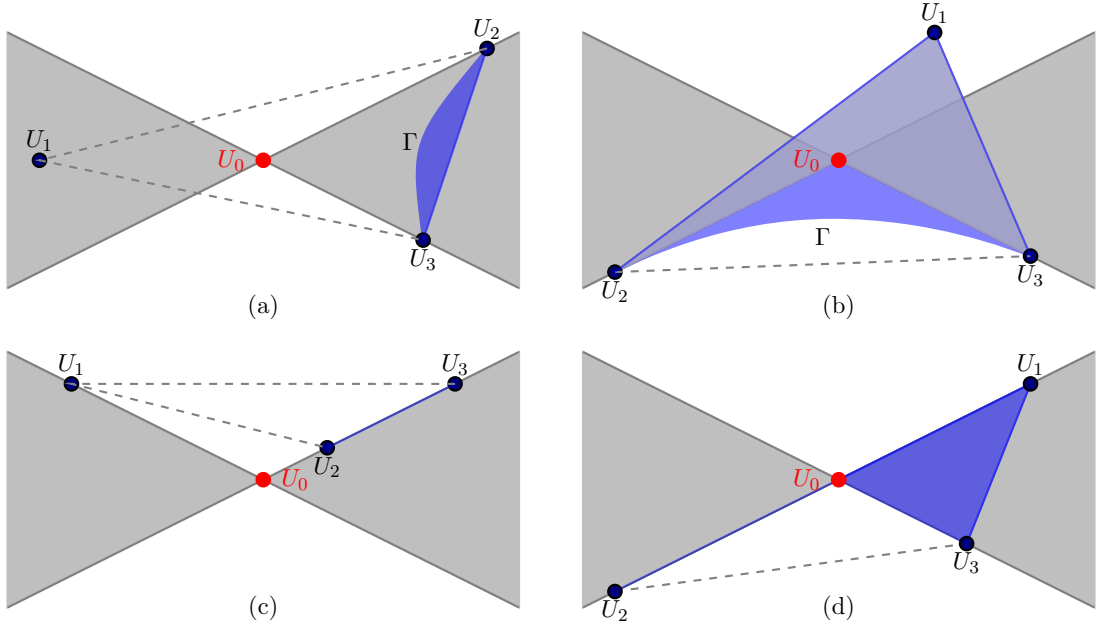
\begin{figure}[h]
\centering
\begin{tikzpicture} [thick,scale=0.85, every node/.style={transform shape}]
\begin{scope}
\filldraw[lightgray] (4,2) -- (8,4) -- (8,0) -- cycle;
\filldraw[lightgray] (4,2) -- (0,4) -- (0,0) -- cycle;
\draw[thick,gray] (0,0) --(8,4);
\draw[thick,gray] (0,4) --(8,0);
\filldraw[blue, opacity=0.5] (7.5,3.75) .. controls (6.3,2.3) and (6.3,2.3) .. (6.5,.75);
\node at (6.3,2.3) {$\Gamma$};
\filldraw[red] (4,2) circle [radius=3pt];
\node[red] at (3.5,2) {$U_0$};
\filldraw[black] (7.5,3.75) circle [radius=3pt] node [anchor=south]{$U_2$};
\filldraw[black] (6.5,.75) circle [radius=3pt]node [anchor=north]{$U_3$};
\filldraw[black] (.5,2) circle [radius=3pt]node [anchor=south]{$U_1$};

\filldraw[blue, opacity=0.5] (7.5,3.75) circle [radius=2pt];
\filldraw[blue, opacity=0.5] (6.5,.75) circle [radius=2pt];
\filldraw[blue, opacity=0.5] (.5,2) circle [radius=2pt];

\draw[blue, opacity=0.5, thick] (7.5,3.75) -- (6.5,.75);
\draw[gray, dashed, thick] (7.5,3.75) -- (.5,2) --(6.5,.75);
\node at (4,-0.3) {(a)};
\end{scope}

\begin{scope}[shift={(9,0)}]


\filldraw[lightgray] (4,2) -- (8,4) -- (8,0) -- cycle;
\filldraw[lightgray] (4,2) -- (0,4) -- (0,0) -- cycle;

\draw[thick,gray] (0,0) --(8,4);
\draw[thick,gray] (0,4) --(8,0);
\filldraw[blue, opacity=0.5] (0.5,0.25) parabola bend (3.9,1.1) (7,0.5) -- (5.5,4) -- cycle;
\node at (3.75,0.8) {$\Gamma$};
\filldraw[lightgray, opacity=0.7] (5.5,4) -- (0.5,.25) -- (4,2) -- (7,0.5) -- cycle;

\filldraw[red] (4,2) circle [radius=3pt];
\node[red] at (3.5,2) {$U_0$};
\filldraw[black] (0.5,.25) circle [radius=3pt] node [anchor=north]{$U_2$};
\filldraw[black] (7,0.5) circle [radius=3pt]node [anchor=north]{$U_3$};
\filldraw[black] (5.5,4) circle [radius=3pt]node [anchor=south]{$U_1$};

\filldraw[blue, opacity=0.5] (0.5,0.25) circle [radius=2pt];
\filldraw[blue, opacity=0.5] (7,0.5) circle [radius=2pt];
\filldraw[blue, opacity=0.5] (5.5,4) circle [radius=2pt];

\draw[gray, dashed, thick] (.5,.25) -- (7,0.5);
\draw[blue, opacity=0.5, thick] (.5,.25) -- (5.5,4) --(7,0.5);
\node at (4,-0.3) {(b)};
\end{scope}

\begin{scope}[shift={(0,-5)}]
\filldraw[lightgray] (4,2) -- (8,4) -- (8,0) -- cycle;
\filldraw[lightgray] (4,2) -- (0,4) -- (0,0) -- cycle;
\draw[thick,gray] (0,0) --(8,4);
\draw[thick,gray] (0,4) --(8,0);
\filldraw[red] (4,2) circle [radius=3pt];
\node[red] at (4.5,2) {$U_0$};
\filldraw[black] (7,3.5) circle [radius=3pt] node [anchor=south]{$U_3$};
\filldraw[black] (5,2.5) circle [radius=3pt]node [anchor=north]{$U_2$};
\filldraw[black] (1,3.5) circle [radius=3pt]node [anchor=south]{$U_1$};

\filldraw[blue, opacity=0.5] (7,3.5) circle [radius=2pt];
\filldraw[blue, opacity=0.5] (5,2.5) circle [radius=2pt];
\filldraw[blue, opacity=0.5] (1,3.5) circle [radius=2pt];

\draw[blue, opacity=0.5, thick] (7,3.5) -- (5,2.5);
\draw[gray, dashed, thick] (7,3.5) -- (1,3.5) --(5,2.5);
\node at (4,0) {(c)};
\end{scope}

\begin{scope}[shift={(9,-5)}]
\filldraw[lightgray] (4,2) -- (8,4) -- (8,0) -- cycle;
\filldraw[lightgray] (4,2) -- (0,4) -- (0,0) -- cycle;
\draw[thick,gray] (0,0) --(8,4);
\draw[thick,gray] (0,4) --(8,0);
\filldraw[blue, opacity=0.5] (4,2) -- (7,3.5) -- (6,1) -- cycle;

\draw[gray, dashed, thick] (.5,.25) -- (6,1);
\draw[blue, opacity=0.5, thick] (.5,.25) -- (7,3.5) --(6,1);
\node at (4,0) {(d)};

\filldraw[red] (4,2) circle [radius=3pt];
\node[red] at (3.5,2) {$U_0$};
\filldraw[black] (0.5,.25) circle [radius=3pt] node [anchor=north]{$U_2$};
\filldraw[black] (6,1) circle [radius=3pt]node [anchor=north]{$U_3$};
\filldraw[black] (7,3.5) circle [radius=3pt]node [anchor=south]{$U_1$};

\filldraw[blue, opacity=0.5] (0.5,0.25) circle [radius=2pt];
\filldraw[blue, opacity=0.5] (6,1) circle [radius=2pt];
\filldraw[blue, opacity=0.5] (7,3.5) circle [radius=2pt];

\end{scope}

\end{tikzpicture}
\caption{
Four different three-well configurations related \cref{T3} and \cref{T2}. Black and red dots represent the wells and $U_0$, respectively. Solid lines represent compatibility between wells and dashed lines incompatibility.    The gray planar cone represents $\Aff(\mathcal{U})\cap \mathcal{C}_{U_0}$, namely the set of incompatible matrices with $U_0$ in $\Aff(\mathcal{U})$.  Gray lines are the rank-one lines across $U_0$. In figure (a) and (b), the blue regions are the outer bounds of $Q^e({\mathcal U})$ in  \cref{T2}. Figures (c) and (d) represent two admissible configurations in \cref{T3}, and blue regions are $Q^e({\mathcal U})=L^e({\mathcal U})$.
}
\label{fig:1}
\end{figure}

We present  \cref{fig:1}{\color{DarkOrchid}.(b)} and \cref{fig:1}{\color{DarkOrchid}.(d)} to display the analogous configurations when only one pair of wells in ${\mathcal U}$  is incompatible. 
In the configuration presented in  \cref{fig:1}{\color{DarkOrchid}.(b)}, $L^e(\mathcal{U})$ is the wedge-like region bounded by the polygon with vertices $U_1,\, U_2,\, U_0$ and $U_3$ by \cref{P:lamconv}, and  the well $U_0\in C(\mathcal{U})$ is rank-one compatible with $U_2$ and $U_3$ simultaneously.  \cref{T2}'s outer bound correspond to the region bounded by the curve $\Gamma$ and the two segments joining $U_1$ with $U_2$ and $U_3$, respectively. The configuration presented in \cref{fig:1}{\color{DarkOrchid}.(d)} has a rank-one compatible pair, so $Q^e(\mathcal{U})$ equals $L^e(\mathcal{U})$ due to \cref{T3}, and both sets equal the blue flag-like region.
\medskip

The organization of the paper is as follows. In \cref{sec:3}, we briefly describe the ideas for the proofs and give some further comments. 
 In \cref{sec:4}, we introduce the notion of incompatible cone and study some 
of its properties. The proof of \cref{P:lamconv}
is presented in \cref{sec:5}. In \cref{sec:6}, we present the symmetric polyconvex conjugate and symmetric biconjugate functions, and provide its application to our particular setting in \cref{sec:7}. 
The proofs of Theorems~\ref{T3} and \ref{T2} are presented \cref{sec:8}. Finally, in \cref{sec:9} we discuss the optimality of the outer bounds in \cref{T2} with respect to symmetric quadratic polyconvex functions.  

\section{Ideas of proofs and further comments}\label{sec:3}
In this paper, we construct quasiconvex hulls' inner and outer bound to $Q^e(\mathcal{U})$. The inner bound is the symmetric lamination convex hull $L^e({\mathcal U})$ that we derive explicitly in \cref{P:lamconv}.  For the outer bound, let  $P^e(\mathcal{U})$ be the symmetric polyconvex hull of ${\mathcal U}$. If $f$ is a non-negative symmetric polyconvex function such that $\mathcal{U}$ belongs to  $\Ker f:=\{M\in \msym\, |\, f(M)=0\}$, then $Q^e(\mathcal{U}) \subset P^e(\mathcal{U}) \subset \Ker f\cap C(\mathcal{U})$. 
The later inclusion is clear from the characterization of $P^e(\mathcal{U})$ in terms of cosets, i.e., \cref{eq:cosets}.  We do not need to construct the polyconvex function $f$ explicitly, but only its kernel. Finally, we construct the required matrix C in terms of the elements in U. With this matrix $C$, $\Ker f^{pp}$ equals the optimal outer bound when only quadratic polyconvex functions are used to bound $Q^e(\mathcal{U})$.

There are not many examples of polyconvex functions, and we need them to vanish on $\mathcal{U}$. Therefore, we construct the parametric family as follows. First, we consider a function of the form  
\[
f_C(M) = \chi_{\bar{B}}(e(M))\left(|\det e(M)| + |\pro{C}{M}|\right), \quad \bar{B} = \msym \setminus L^e(\mathcal{U}),
\] 
where $\det C<0$, and  $\chi_{\bar{B}}:\msym\rightarrow \{0,1\}$ is the indicator function of the set $\bar{B}$. Notice that by construction,  the functions $f_C$ vanish on  $L^e(\mathcal{U})$ (hence on $\mathcal{U}$), and they are quadratic outside this set. 
Second, we compute the symmetric polyconvex envelope of $f_C$ through its symmetric polyconvex  biconjugate function $f_C^{PP}$. Then, we estimate   $\Ker f_C^{PP}$  in  \cref{eqboundary}.

In our results, we describe the geometry of the laminar convex and quasiconvex hulls in terms of the matrix $U_0$. This matrix is defined as the intersection of two rank-one lines passing through two  wells in ${\mathcal U}$ (see \cref{fig:1}). These lines exist since the cone ${\mathcal C}_{0}$ centered at any element of ${\mathcal U}$ intersects the plane defined by $\Aff({\mathcal U})$.  If we denote the normal matrix to $\Aff({\mathcal U})$ by $Q$, we show that the later condition is equivalent to $\det Q<0$ (see \cref{rmk:det}).  
This is related to the condition $\det C < 0$ in the construction of $f_C$. Therefore, in our arguments, the geometric relationship between the incompatible cone and $\Aff({\mathcal U})$ is important. 

In the last section, we show that the outer bound in \cref{T2} is optimal for quadratic symmetric polyconvex functions. Indeed, we notice that this bound  depends on the quadratic polyconvex function $\hbar$, and the proof of the optimality for quadratic polyconvex functions heavily depends
on the recent characterization of symmetric polyconvex functions in two dimensions given by Boussaid {\em et. al} \cite{Anja19} (proposition 4.5 p.435), and it can not be extended beyond quadratic polyconvexity without an analogous result for non-quadratic functions.  
We believe that better estimates can be obtained by considering extensions of the form 
\begin{equation*}\label{fCp}
f_C(M) = \chi_{\bar{B}}(e(M))\left(|\det e(M)|^q + |\pro{C}{M}|^p\right), 
\end{equation*}
for some  positive exponents $p$ and $q$. The computation in these cases become  cumbersome, and we do not pursue them further in this paper. Another consequence of the optimality in the outer bound for quadratic functions is that, if the outer bounds of $Q^e({\mathcal U})$ in \cref{T2} are not sharp, then the quasiconvex envelope of an energy density $f$ (as in \cref{freeEC}), is not quadratic. 

Regarding possible extensions of the presented results, we study a generalization of \cref{T3}  to four and more co-planar wells in two dimensions in  \cite{camo2}. In particular, we are interested in the sufficient conditions on rank-one compatible pairs  to show that $Q^e(\mathcal U)=L^e(\mathcal U)$. 

Finally, we make a few comments on the possible extension of our results to the three-dimensional case. The first step is to identify the corresponding geometry. We need to generalize the definition of the cone ${\mathcal C}_0$ and the normal $Q$ to $\Aff({\mathcal U})$.  To generalize the cone ${\mathcal C}_0$, we notice that two matrices 
$M_1$ and $M_2$ in $\re^{3\times 3}_{sym}$ represent compatible linear strains if the eigenvalues of $M_1-M_2$ are of the form $\lambda_1 <\lambda_2 =0 < \lambda_3$. Now, the cone's boundary is a
``surface'' that represents the rank-one compatible directions. Therefore, by the isomorphism between $\re^{3\times 3}_{sym}$ and $\re^6$, this boundary should be an object of four dimensions in $\re^6$.
Consequently, the tangent space to the cone's boundary is of dimension four, and we conjecture that we need at most five wells to set up the analogous geometry. The challenge, in this case, is to define the analogous to $Q$ and the generalization of condition $\det Q <0$ that identifies the intersection of the incompatible cone and the affine space. In some cases, we believe that it is possible that for a particular set of wells in three dimensions, the cone geometry reduces to its two-dimensional counterpart, and we recover results similar to the ones presented in this paper. As a piece of evidence in this direction, we recall the results in \cite{Anja2015} quoted before.

\section{On the geometry of the incompatible cone}\label{sec:4}

We begin with an equivalence result on compatibility. 
\begin{proposition}\label{qconnectness}
Let $M_1,M_2\in \m$. 
The following statements are equivalent:
\begin{enumerate}[label=(\alph*)]
\setlength\itemsep{.2em}
\item\label{qcon1}  There exist $v\in \re^2$, $a \in S^1$ and $W\in \mskew$ such that $M_1-M_2 + W =a^\perp\otimes v.$ 
\item\label{qcon2} There exist $v\in \re^2$ and $a \in S^1$ such that $e(M_1) = e(M_2) + v\odot a^\perp$. 
\item\label{qcon3}  There exists $a \in S^1$ such that  $\pro{e(M_1)-e(M_2)}{a\otimes a} = 0$. 
\item\label{qcon4}   $\det (e(M_1)-e(M_2))\le 0$.
\end{enumerate}
\end{proposition}

\begin{proof} 
That \ref{qcon1} implies \ref{qcon2} follows directly from the symmetric parts of the equation in the statement \ref{qcon1}. 
Since $a^\perp$ is orthogonal to $a$, we easily conclude that  \ref{qcon2} implies \ref{qcon3}. To show that  \ref{qcon3} implies  \ref{qcon4}, we proceed by contradiction and  assume $\det (e(M_1)-e(M_2))>0,$ meanwhile \ref{qcon3} holds. Then it follows that $a\mapsto (e(M_1)-e(M_2))a\cdot a$ is either a positive or negative defined quadratic form. Thus, 
the unique solution to $(e(M_1)-e(M_2))a \cdot a=0$ is  $a=0$, a contradiction to $a\in S^1$ by hypothesis.

Now,  we show that \ref{qcon4} implies \ref{qcon1}. 
Since $\det M = \det(e(M))+\det(w(M))$ for every $M\in\m$ and \ref{qcon4} is assumed, we obtain 
 \[
\det(M_1-M_2 +\mu R) = \det(e(M_1)-e(M_2)) +\left(\frac{1}{2}\pro{M_1-M_2}{R}+\mu\right)^2.
\] 
Thus, the equation $\det(M_1-M_2 +\mu R)=0$ is satisfied by some   
$\mu^*\in \re$. 
So,  the  statement \ref{qcon1} follows, 
and the proof is finished.  
 \end{proof}

 The following lemma states that if $M\in \msym$ has non-positive determinant, then it is the symmetric part of a tensor product, this result is used in the proof of Lemma \ref{planeQV}.

\begin{lemma}\label{det<}
Let $Q\in \msym$ be such that its determinant is non-positive. There exist $a, \,n\in S^1$ and $\nu\in \re$ 
such that $Q = \nu \sten{a}{n}$ and $\det Q = -\nu^2(a\cdot n^\perp)^2/4$. Moreover, $v_+=a+n$, and 
$v_-=a-n$ are eigenvectors of $Q$ with eigenvalues $\lambda_+ = \nu[a\cdot e+1]/2$, and 
$\lambda_- = \nu[a\cdot e-1]/2,$ respectively.
\end{lemma}
\begin{proof}
The existence of $a, \,n\in S^1$, and $\nu\in \re$ such that $Q = \nu \sten{a}{n}$ follows by setting 
$M_1=Q$, $M_2=0$ into the equivalence between \ref{qcon2} and \ref{qcon4} in \cref{qconnectness}.  
The statement concerning the eigenvalues and eigenvectors follows by a direct computation
\[
Qv_\pm = \frac{\nu}{2}(a\otimes n + n\otimes a)(a\pm n) =\frac{\nu}{2}[(a\cdot n)(a\pm n)+(n\pm a)] = \frac{\nu}{2}[a\cdot n \pm 1](a\pm n).
\]
Thus, $\det Q = \lambda_+ \lambda_- = -\nu^2[1-(a\cdot n)^2]/4 = -\nu^2\left (a\cdot n^\perp \right)^2/4$. 
\end{proof}
The next lemma characterizes the set of symmetric  incompatible matrices as elements in a solid cone. In particular, the set of symmetric matrices that are rank-one compatible is identified with the boundary of this cone.
\begin{lemma}\label{conenonconnection}
 Let $M\in \msym$, then  $\det M = 0$ if and only if $|\pro{M}{\Id}| = \|M\|$. Moreover, $\det M < 0$ if and only if $|\pro{M}{\Id}| < \|M\|$.
\end{lemma}
\begin{proof} 
Let  $M\in\msym$ be given by 
 \[
 M = \begin{pmatrix}
 M_{11} & M_{12} \\ M_{12} & M_{22}
 \end{pmatrix}
 \]
with $ M_{11}, M_{12}, M_{22}\in \re$, and consider the square of its Frobenius norm 
\[
\begin{array}{rl}
\|M\|^2 
 = & M_{11}^2 + M_{22}^2 + 2M_{12}^2
 =  (M_{11} + M_{22})^2 + 2(M_{12}^2 -M_{11}M_{22}). 
\end{array}
\]
Hence, 
$\|M\|^2 =  \pro{M}{\Id}^2 - 2\det M$,
and the proof follows. 
\end{proof}

\begin{figure}[h]
\centering
\includegraphics[scale=.8]{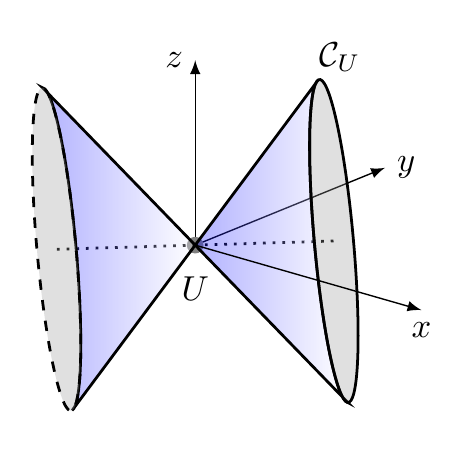}
\caption{The incompatible cone $\mathcal C$ represented in $\re^3$ by the isomorphism \eqref{isomorphism}}
\label{fig:coneR3}
\end{figure}

Because of \cref{qconnectness} and \cref{conenonconnection}, the set of incompatible linear strains with $U\in \msym$ is given by the interior of the cone 
\[ 
\begin{split}
\mathcal{C}_U :=& \left\{V\in\msym \mbox{ such that } \  \|V-U\| < |\pro{V-U}{\Id}| \right \},\\
=&\left\{V\in\msym \mbox{ such that } \  \frac{1}{\sqrt{2}} < \left|\pro{\frac{V-U}{\|V-U\|}}{\frac{1}{\sqrt{2}}\Id}\right| \right \}
\end{split}
\] 
Let $e^1$, $e^2$, and $e^3$ stand for the canonical vectors in $\mathbb{R}^3$. Then $e^1\otimes e^1$ and $e^2\otimes e^2$ belong to $\partial \mathcal{C}_0$, and the cone's aperture angle, $2\theta$, is $\pi/2$ since $e^1\otimes e^1$, $e^2\otimes e^2$ and $\Id$ are coplanar, see \cref{fig:coneR3}.

\begin{remark}\label{det:prop}
Let ${\mathcal S} :\m \rightarrow \m$ stand for the linear map $M \mapsto RMR^T$, with $R$ as in \eqref{Rmatrix}. The function ${\mathcal S}$ maps the  symmetric  (skew-symmetric) subspace into itself. Moreover, a simple calculation shows that $adj(M) = \mathcal{S}M$ for every $M\in \m$, and we get the following identities: 
\begin{enumerate}[label=(\alph*)]
\setlength\itemsep{.2em}\item\label{det:prop1} $ 2\det M = \pro{\mathcal{S}M}{M}.$ 
\item\label{det:prop2} $\det M = \det e(M) + \det w(M)= \det e(M) + w^2_M.$
\item\label{det:prop3} $\det (N+M) = \det M + \det M + \pro{\mathcal{S}M}{N},$ \quad $N\in \m.$
\item\label{det:prop4} If $M = \pv{\xi}{\eta}$ for some $a,\,n\in S^1$ and $\xi,\,\eta\in \re$, then $\det M = \xi \eta |a\times n|^2$.
\end{enumerate}
\end{remark}
 
The next lemma describes a coordinate system defined by the intersection between the boundary of the incompatible cone and $\Aff({\mathcal U})$.
 This system of coordinates is convenient and simplifies many of the computations in our proofs.
\begin{lemma}\label{planeQV}
Let $V, Q\in \msym$ such that $\det Q<0$ and denote by $\Pi_Q(V)$ the affine set of codimension one that contains $V$ and is normal to $Q$,  namely 
\[\Pi_Q(V) = \{U\in \msym\,|\,\pro{Q}{U-V} =0\}.\] Then, there exist two nonparallel vectors $a,n\in S^1$, such that 
\begin{equation}\label{planeQV-1}
\Pi_Q(V) = \{V + \xi\, \ten{a^\perp}{} + \eta\,\ten{n^\perp}{} \,| \, (\xi,\eta)\in \re^2\}.
\end{equation}
Moreover, $U\in \Pi_Q(V)$ and $V$ represent compatible linear strains if and only if $U - V = \xi \,\ten{a^\perp}{} + \eta\,\ten{n^\perp}{}$ with  $\xi\,\eta\leq 0$.
\end{lemma}

\begin{proof}
By Lemma \ref{det<}, there exist $a,n\in S^1,$ and $\nu\in \re \setminus\{0\}$  such that $Q=\nu a\odot n$. 
The vectors $a$ and $n$ are linearly independent since $\det Q<0$.
Moreover,  by  \cref{qconnectness}~\ref{qcon2} and \ref{qcon3}, and the symmetric role between 
vectors $a$ and $n$, we have  that $\ten{a^\perp}{}$ and $\ten{n^\perp}{}$  are two  linearly independent rank-one matrices such that
\begin{equation*}
\pro{Q}{\ten{a^\perp}{}} = \pro{Q}{\ten{n^\perp}{}} = 0. 
\end{equation*}
Hence, $\{\ten{a^\perp}{},\,\ten{n^\perp}{}\}$ is a basis for the two-dimensional subspace $\Pi_Q(V)-V$ and \eqref{planeQV-1} follows. 
Now, let $W\in \Pi_Q(V)$ and $(\xi,\eta)\in \re^2$  such that $W-V = \xi\, \ten{a^\perp}{} + \eta\,\ten{n^\perp}{}$. Thus $\det(W-V) = \xi\eta |a\times n|^2$, and 
the second part of the statement follows straight forward from \cref{qcon4} in \cref{qconnectness}. 
\end{proof}

 \begin{remark}\label{rmk:det}
It is well known, see \cite{BFJK,BDa} pages p.191 and p.825 respectively, that if $f:\re^{m\times m}\rightarrow \re$ is quadratic and $f(a\otimes n)\geq 0$ for every $a,n\in\re^m$ then $f$ is rank-one convex. Thus, $-\det(\cdot):\re^{2\times 2}_{sym}\rightarrow \re$ is a symmetric rank-one convex function because \cref{det:prop1} in \cref{det:prop} implies that it is quadratic and the equivalence between items \ref{qcon2} and \ref{qcon4} in \cref{qconnectness} implies that 
$-\det(e(a\otimes n))\geq 0$ for every $a,n\in\re^2$. 
\end{remark}

\begin{lemma}\label{lem:linedot}
Let $U,V,W\in \re^{2\times 2}_{sym}$ such that $V,W$ are compatible and  $U$ is incompatible with both of them. Then $U$ is incompatible with any point in $C(\{U,V,W\})\setminus \{U\}$.  
 \end{lemma}
\begin{proof}
Let $M\in C(\{U,V,W\})\setminus \{U\}$. 
Then, $M=\lambda_1U + \lambda_2V + \lambda_3W$ for some $\lambda_1\in[0,1)$ and $\lambda_2,\lambda_3\in [0,1]$,  such that $\lambda_1+\lambda_2+\lambda_3 =1$.
Since $\lambda_1\neq 1$, we have that $\lambda_2 + \lambda_3>0$ and 
$$
M-U =(\lambda_2+\lambda_3)\left(s(V-U) + (1-s)(W-U)\right),
$$
where $s=\lambda_2/(\lambda_2+\lambda_3)$, and $(1-s)= \lambda_3/(\lambda_2+\lambda_3)$.
Hence, $\det(M-U) = (\lambda_2+\lambda_3)^2 \det(s(V-U) + (1-s)(W-U))$. 
Now, since $V-U$ and $W-U$ are compatible, and $-\det(e(\cdot))$ is rank-one convex, we have that 
\[
\det(s(V-U) +(1-s)(W-U)) \geq s\det(V-U)+(1-s)\det(W-U)>0.
\]
Therefore, $U$ and $M$ are incompatible and the proof is complete.
\end{proof}

In Proposition~\ref{P:lamconv} and Theorem~\ref{T2}, we consider two possible configurations of wells depending on the number of compatible pairs among its elements. We use the following definition to keep our presentation compact.  
\todoin{
\begin{definition}
Let $\mathcal{U}$ represent a three-well set such that $\Aff(\mathcal{U})$ has codimension one. We say that $\mathcal{U}$ is {\em type one} if up to a relabeling $\mathcal{U}=\{U_1,U_2,U_3\}$, where  
\begin{equation}\label{eq:1class}
\det(U_1-U_2)> 0, \ \det(U_1-U_3)> 0, \mbox{ and }\det(U_2-U_3) \leq 0,
\end{equation}
and $\mathcal{U}$ is {\em type two} if up to a relabeling $\mathcal{U}=\{U_1,U_2,U_3\}$, where  
\begin{equation}\label{eq:2class}
\det(U_1-U_2)\leq 0,\ \det(U_1-U_3)\leq 0, \mbox{ and }\det(U_2-U_3)>0.
\end{equation}
\end{definition}
\begin{remark}\label{rem:det<}
We claim that if $\mathcal{U}$ is either type one or type two  and $Q$ is normal to $\Aff(\mathcal{U})$, then $\det Q<0$. We argue by contradiction. Indeed, if $\det Q\geq 0$, then $\Aff(\mathcal{U})\cap \mathcal{C}_{U_1}= \{U_1\}$ and  $U_1$ is compatible with $U_2$ and $U_3$,  see paragraph after \cref{conenonconnection}. We conclude that  $\mathcal{U}$ is pairwise compatible by arguing the same for the two remaining wells, but this is a contradiction to \cref{eq:1class} and \cref{eq:2class}.
\end{remark}
}

In the following lemma, we apply the coordinate system in \cref{planeQV}  to the geometry of our three-well problem.  
\begin{lemma}\label{origin}
Let $\mathcal{U}\subset \msym$ represent a three-well set where $\mathcal{U}$ is either type one or type two. Then, there exist $U_2,U_3\in \mathcal{U}$ and  $U_0\in C(\mathcal{U})$ such that $\det(U_2-U_0)=\det(U_3-U_0)=0$.
\end{lemma}

\begin{proof}
First, by \cref{rem:det<}, $\det Q<0$, and by \cref{det<}, there exist $a,n \in S^1$ linearly independent vectors such that $Q=\nu a\odot n$ for some $\nu\in \re\setminus \{0\}$. Hence, by Lemma \ref{planeQV}, 
\begin{equation}\label{U2U3}
U_2= U_1 + \xi_2 a^\perp\otimes a^\perp + \eta_2 n^\perp\otimes n^\perp, \quad \mbox{and} \quad U_3= U_1 + \xi_3 a^\perp\otimes a^\perp + \eta_3 n^\perp\otimes n^\perp,
\end{equation}
for some $(\xi_2,\eta_2), (\xi_3,\eta_3)\in \re^2$. We notice that
$(\xi_2,\eta_2)$ and $(\xi_3,\eta_3)$ belong to the same quadrant in the $\re^2$. Otherwise, all elements in $\mathcal{U}$ are either compatible or incompatible and this contradicts either \cref{eq:1class} or \cref{eq:2class}.

Second, we claim that 
\begin{equation*}
U_0 = U_1 +\alpha a^\perp\otimes a^\perp + \beta n^\perp\otimes n^\perp, \quad \mbox{with}\quad \begin{cases} \alpha = \mbox{arg} \min \{|\xi|\, |\, \xi\in\{\xi_2, \xi_3\}\},\\
\beta = \mbox{arg} \min \{|\eta|\, |\, \eta\in\{\eta_2, \eta_3\}\},
\end{cases}
\end{equation*}
is the desired matrix. Indeed, without loss of generality we assume ${\mathcal U}$ is type one. Then, either $\alpha=\xi_3$ or $\beta=\eta_3$, but not both. Otherwise, by Remark \ref{det:prop}~\ref{det:prop4} and \cref{U2U3}, we get $\det(U_3-U_2)=(\xi_3-\xi_2)(\eta_3-\eta_2)|a\times n|^2\geq 0$, a contradiction with ${\mathcal U}$ being of type one. 
Hence, $\det(U_3-U_0) = (\xi_3-\alpha)(\eta_3-\beta)|a\times n|^2=0$, and by an analogous argument, 
we also have $\det(U_2-U_0) =0$. Now, since $U_1$ is incompatible with both $U_2$ and $U_3$  (namely, $\xi_2\eta_2>0$ and $\xi_3\eta_3>0$), and $(\xi_2,\eta_2)$ and $(\xi_3,\eta_3)$ belong to the same quadrant, we have that $\det(U_1-U_0)=\alpha\beta|a\times n|^2\geq 0$. We let
{\footnotesize \[
\lambda_1 =\frac{(\alpha-\xi_2)(\beta-\eta_3)-(\alpha-\xi_3)(\beta-\eta_2)}{\xi_2\eta_3-\xi_3\eta_2}, \ \ \lambda_2=\frac{(\alpha-\xi_3)\beta-(\beta_-\eta_3)\alpha}{\xi_2\eta_3-\xi_3\eta_2} , \ \ \mbox{and} \ \ \lambda_3=\frac{(\beta-\eta_2)\alpha-(\alpha-\xi_2)\beta}{\xi_2\eta_3-\xi_3\eta_2}.
\]}
By a straight forward computation, we
get that $U_0=\lambda_1U_1 + \lambda_2U_2 + \lambda_3U_3$ with 
$\lambda_1, \lambda_2, \lambda_3 \in (0,1)$, and $\lambda_1+\lambda_2 + \lambda_3 = 1$. Hence, $U_0\in C(\mathcal{U})$ as claimed. The case where ${\mathcal U}$ is type two can be treated analogously,  and the result follows. 
\end{proof}

\begin{corollary}\label{cor:origin}
Let $\mathcal{U}\subset \msym$ represent a three-well set where  $\mathcal{U}$ is either type one or type two. Then, there exist $U_0\in C(\mathcal{U})$, $a,n\in S^1$, and $\xi,\ \eta,\ \gamma,\ \zeta \in \re$ such that  
\begin{equation}\label{repUs}
U_1 = U_0 +\pv{\xi}{\eta}, \quad U_2 = U_0+\pva{\gamma}, \quad \mbox{and} \quad U_3 = U_0 + \pvn{\zeta},
\end{equation}
where 
\begin{equation}\label{coefU}
\begin{cases}
\eta>0\geq \zeta,\ \ \xi>0\geq \gamma \ \ \mbox{or} \ \ \zeta \geq 0>\eta,\ \ \gamma\geq 0>\xi 
&\mbox{if }\mathcal{U} \mbox{ is type one,}\\
\eta\geq0> \zeta,\ \ \gamma>0\geq \xi \ \ \mbox{or} \ \ \zeta > 0\geq \eta,\ \ \xi\geq 0>\gamma 
&\mbox{if }\mathcal{U} \mbox{ is type two.} 
\end{cases}
\end{equation}

\end{corollary}
\begin{proof}
By \cref{origin}, there exists $U_0\in C(\mathcal{U})$ that is rank-one compatible with $U_2$ and $U_3$. 
From, \cref{rem:det<}, the normal $Q$ to $\Aff(\mathcal{U})$ is such that $\det Q<0$. Then, \eqref{repUs} follows from 
\cref{det<}, \cref{planeQV}, and  that $\Pi_Q(U_0)=\Aff(\mathcal{U})$ has codimension one. 

Next, we prove \cref{coefU}. Without loss of generality we assume that ${\mathcal U}$ is of type two. Hence, $\det(U_2-U_3)>0,$ and $\det(U_j-U_1)\leq 0$ for $j=0,2,3$.
Therefore, by \eqref{repUs} and \cref{det:prop}~\cref{det:prop4}, we get
\[
-\gamma\,\zeta>0,\quad \xi\,\eta\leq 0,\quad (\xi-\gamma)\eta\leq 0,\quad \mbox{and}\quad (\eta-\zeta)\xi\leq 0.
\]
Now, the conditions in \eqref{coefU} follows since  $U_0\in C(\mathcal{U})$. The case where $\mathcal{U}$ is type one can be proved by a similar argument. 
\end{proof}
\section{The symmetric lamination convex hull of \texorpdfstring{$\mathcal{U}$}{}}\label{sec:5}

In this section, we determine the symmetric lamination convex hull when $\mathcal{U}$ is either type one or type two, {\em i.e.} there is one out of the three wells in $\mathcal{U}$ that is compatible with one of the two remaining wells and incompatible with the other.

\begin{proof}{\em of \cref{P:lamconv}:}
We prove the item \ref{lhull1}. 
By hypothesis, $U_2$ is compatible with $U_3$, and $U_1$ is incompatible with both of them. 
Then, $L^{e,1}(\mathcal{U}) =\{U_1\} \cup C(\{U_2,U_3\})$. By \cref{lem:linedot} there is not $U\in C(\mathcal{U})$ compatible with $U_1$.
Therefore, no laminate of degree two is admissible and $L^{e,2}(\mathcal{U})=L^{e,1}(\mathcal{U})$. 
By the same argument, $L^{e,n}(\mathcal{U})=L^{e,1}(\mathcal{U})$ for every $n\ge 1$, and \cref{lhull1} follows from the definition of $L^{e}(\mathcal{U})$.

Now, we prove the item \ref{lhull2}. In this case  $\mathcal{U}$ is of type two. Then,  \cref{origin} implies the existence of $U_0\in C(\mathcal{U})$ such that
\begin{equation*}
\det(U_2-U_0)=\det(U_3-U_0)=0,
\end{equation*}
 and we assume that the wells in $\mathcal{U}$ are given by \cref{repUs} in \cref{cor:origin}.
Now, by the compatibility relations, $\mathcal{U}$'s first symmetric lamination is $L^{e,1}(\mathcal{U}) = C(\{U_1,U_2\})\cup C(\{U_1,U_3\})$.  For the laminations of degree
greater than one, first, we assume that 
$\det(U_1-U_2)<0$ and $\det(U_1-U_3)<0$ and consider the continuous function 
\[
t\mapsto \det(tU_1+(1-t)U_3-U_0) = t^2\det(U_1-U_3)+t\pro{\mathcal{S}(U_3-U_1)}{U_3-U_0}.
\]
By \cref{det:prop}, \cref{repUs}, and \eqref{coefU}~(b) for the type two case, we have that
\[
\pro{\mathcal{S}(U_3-U_1)}{U_3-U_0}
=  -\xi\zeta |n\times a|^2 >0.
\]
Thus, there exists $t_0>0$ such that $U_{1,3}=t_0U_1+(1-t_0)U_3\in C(\{U_1,U_3\})$ and $\det(U_{1,3}-U_0)=0$ since $\det(U_1-U_3)<0$. Analogously, there exists $U_{1,2}\in C(\{U_1,U_2\})$ such that $\det(U_{1,2}-U_0)=0$.
\begin{figure}
\centering
\begin{tikzpicture} [thick,scale=0.85, every node/.style={transform shape}]
\begin{scope}[shift={(0,0)}]
\filldraw[lightgray] (4,2) -- (8,4) -- (8,0) -- cycle;
\filldraw[lightgray] (4,2) -- (0,4) -- (0,0) -- cycle;
\filldraw[brown] (5.94,2.97) circle [radius=3pt];
\node[brown] at (6.2,3.3) {$U_{1,3}$};
\filldraw[brown] (3.3,2.35) circle [radius=3pt];
\node[brown] at (3.2,2.7) {$U_{1,2}$};

\filldraw[blue] (3.4,.8) circle [radius=3pt];
\node[blue] at (3.4,1.15) {$U$};

\filldraw[blue] (2.4,.3) circle [radius=3pt];
\node[blue] at (2.4,-0.1) {$V_{3}$};

\filldraw[blue] (6.26,2.23) circle [radius=3pt];
\node[blue] at (6.4,2.56) {$V_{1,3}$};
\draw[dotted, gray] (2.4,.3) -- (6.26,2.23);

\filldraw[blue] (4.2,.4) circle [radius=3pt];
\node[blue] at (4.3,0.05) {$V_{2}$};

\filldraw[blue] (2.1,1.45) circle [radius=3pt];
\node[blue] at (2.1,1.8) {$V_{1,2}$};
\draw[dotted, gray] (4.2,.4) -- (2.1,1.45);

\filldraw[red] (4,2) circle [radius=3pt];
\node[red] at (3.5,2) {$U_0$};
\filldraw[black] (0.5,.25) circle [radius=3pt] node [anchor=north]{$U_2$};
\filldraw[black] (7,0.5) circle [radius=3pt]node [anchor=north]{$U_3$};
\filldraw[black] (5.5,4) circle [radius=3pt]node [anchor=south]{$U_1$};
\draw[gray, dashed, thick] (.5,.25) -- (7,0.5);
\draw[gray, thick] (.5,.25) -- (5.5,4) --(7,0.5);
\end{scope}
\end{tikzpicture}
\caption{Auxiliary wells used in the proof of \cref{P:lamconv}. This figure shows a three-well set with strict inequalities in \cref{lhull2}. The incompatible cone of $U_0$ restricted to  $\Aff(\mathcal{U})$ is represented by the shaded region. }
\label{fig:points}
\end{figure}
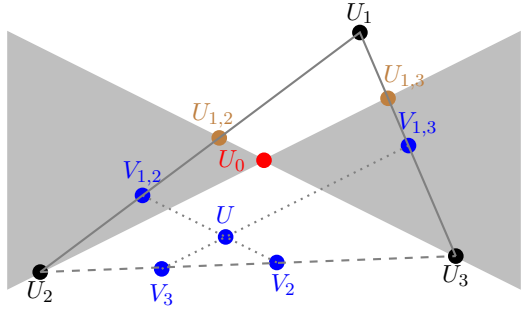
Now, we define two triangular regions with vertices $K_2=\{U_1,U_2,U_{1,3}\}$ and $K_3=\{U_1,U_3,U_{1,2}\}$ (see \cref{fig:points}), and 
divide the proof in three steps:
\begin{enumerate}[label=(\roman*)]
\item\label{step1}  
We claim that  $C(K_2),C(K_3)\subset L^{e,2}(\mathcal{U})$. 
Indeed, if  we assume that $M$ belongs to the relative boundary of $C(K_2)$,  then either $M\in C(\{U_1,U_2\})$ or $M\in C(\{U_2,U_{1,3}\})\cup C(\{U_1,U_{1,3}\}) $. 
Since $U_{1,3}\in L^{e,1}(\mathcal{U})$ and $C(\{U_1,U_2\}) \subset L^{e,1}(\mathcal{U})$ then $M\in L^{e,2}(\mathcal{U})$. 
Now, we assume that $M$ belongs to the relative  interior of $C(K_2)$. 
By construction (see \cref{conenonconnection}), any two matrices  on the line $\ell$ joining $U_2$ and $U_{1,3}$, or on any other parallel line, are rank-one compatible.  Then, there exists a line $\ell_M$ parallel to $\ell$ passing through $M$, such that $\ell_M$ intersects $C(\{U_1,U_2\})$ and $C(\{U_1,U_{1,3}\})$ at some points $V$ and $W$ respectively. Thus, $M\in C(\{V,W\})$, $\det(V-W) =0$, and $V,W\in L^{e,1}(\mathcal{U})$.  Hence, $M\in L^{e,2}(\mathcal{U})$.
By a similar argument $C(K_3)\subset L^{e,2}(\mathcal{U})$ and the claim follows.

\item\label{step2} Now, we claim that $C(K_2)\cup C(K_3) = L^{e,2}(\mathcal{U})$. By contradiction, we assume there exists $U\in L^{e,2}({\mathcal U})\setminus \left(C(K_2)\cup C(K_3)\right)$ (see \cref{fig:points}). 
The intersection between $\partial {\mathcal C}_U$ and  the relative boundary of $C({\mathcal U})$ 
yields the existence of  $V_2,V_3\in C(\{U_2,U_3\})$, $V_{1,2}\in C(\{U_1,U_2\})$, and  $V_{1,3}\in C(\{U_1,U_3\})$ such that $\det (V_2-U)=\det(V_{1,2}-U)=0$  and $\det (V_3-U)=\det(V_{1,3}-U)=0$. Therefore, each matrix  in  the set $\mathcal{U}$  that is compatible with $U$ belongs to either the set
$C_1= C(\{U,V_{1,3}, U_1, V_{1,2}\})$ or the set $C_2=C(\{U,V_2,V_3\})$. Now, since 
$U\in L^{e,2}(\mathcal{U}),$ there exist two compatible matrices $V,W\in L^{e,1}(\mathcal{U})$ 
such that $U=\lambda V + (1-\lambda) W$ for some $\lambda\in (0,1)$. Notice that 
$C_1$ and $C_2$ are convex sets and  $C_1\cap  C_2=\{U\}$, so $V$ and $W$ cannot belong both to $C_1$ or $C_2$ simultaneously. The contradiction follows from the fact $C_2\cap L^{e,1}(\mathcal{U})= \emptyset$ since $U_2$ and $U_3$ are incompatible.

\item\label{step3} We claim that $L^{e,2}(\mathcal{U})=L^{e,3}(\mathcal{U})$. 
The proof follows by  a contradiction argument similar to the one used before. 
Indeed, if there exists $U\in L^{e,3}(\mathcal{U})\setminus L^{e,2}(\mathcal{U})$, 
then $U$ is a convex combination of two compatible matrices $V,W\in L^{e,2}(\mathcal{U})$. 
Since $C_1\cap C_2 = \{U\}$, $V$ and $W$ cannot belong both to $C_1$ or $C_2$. 
Again the contradiction emerged from the fact that $C_2\cap L^{e,2}(\mathcal{U})= \emptyset$. 
Hence, $L^{e,3}(\mathcal{U})\subset L^{e,2}(\mathcal{U})$.  
The reverse inclusion follows from the definition of $L^{e,3}(\mathcal{U})$.
\end{enumerate}

Second, we assume that there is only one rank-one compatibility among the elements of $\mathcal{U}$ and without loss of generality, let $\det(U_1-U_2)<0$ and $\det(U_1-U_3) = 0$. In this case, we define $K_2=C(\{U_1,\,U_0,\, U_2\})$ and $K_3=C(\{U_1,\,U_3\})$. Repeating the argument in step \ref{step1}, we easily get that $K_2\subset L^{e,2}(\mathcal{U})$. Moreover, by definition $K_3\subset L^{e,1}(\mathcal{U}) \subset L^{e,2}(\mathcal{U})$, we obtain the claim in step \ref{step1}. 
The steps \ref{step2} and \ref{step3} are exactly the same.

Third, we assume that there are two rank-one compatible pairs in $\mathcal{U}$, so without loss of generality, $\det(U_2-U_1) =\det(U_3-U_1)=0$. For this case, we choose $K_2=C(\{U_1,\, U_2\})$ and $K_3=C(\{U_1,\,U_3\})$. Under these assumptions, $K_1\cup K_2=C(K_1)\cup C(K_2) = L^{e,1}(\mathcal{U})$. Hence, by step \ref{step2}, $L^{e,1}(\mathcal{U}) = L^{e,2}(\mathcal{U})$.  

Finally, due to 
$L^{e,3}(\mathcal{U})=L^{e,2}(\mathcal{U})$ we get that  $L^{e,n}(\mathcal{U})=L^{e,2}(\mathcal{U})$ for every $n\ge 2$ in all cases,  and \cref{lhull2} follows from the definition of $L^{e}(\mathcal{U})$.
 \end{proof}

\section{Results about polyconvex conjugate and biconjugate functions }\label{sec:6}

In this section we presents a method to construct non-negative polyconvex functions through  the polyconvex conjugate and biconjugate functions $f^p(\xi^*):\re^{m_n} \rightarrow \re \cup \{\pm \infty\}$ and  $f^{pp}(\xi):\re^{n\times n} \rightarrow \re \cup \{\pm \infty\}$ of a given function $f:\re^{n\times n}\rightarrow \re$, respectively. These functions are defined as 
\[
f^p(\xi^*):=\sup\{(T(\eta),\xi^*)-f(\eta)\, | \, \eta\in \re^{n\times n} \}, \quad \mbox{and} \quad  f^{pp}(\xi):=(f^p)^*(T(\xi)), \]
where $T:\re^{n\times n}\rightarrow \re^{m_n}$  is given by \eqref{Tadj}, and $(f^p)^*:\re^{m_n}\rightarrow \re\cup\{\pm\infty\}$ is the dual convex of $f^p$, which is given by $ (f^p)^*(\eta):=\sup\{(\eta,\xi^*)-f^p(\xi^*) \, | \, \xi^* \in \re^{m_n}\}$. 
The following theorem was proved by Kohn \& Strang \cite{KS1} (lemma 3.4 p, 133) also, see theorem 6.6 \cite{BDa} p.268.


\begin{theorem}\label{KS}
Let $f:\re^{n\times n}\rightarrow \re$ and $f^{pp}(M):=(f^p)^*\circ T(M)$. If $g:\re^{n\times n}\rightarrow \re$ is a polyconvex function such that $f(\xi)\geq g(\xi)$ for every $\xi\in \re^{n\times n}$, then $f^{pp}=Pf$, the polyconvex envelope of $f$. 

\end{theorem}

In the two-dimensional framework, $f:\m\to\re$ and, the above expressions become
\begin{equation}\label{fconjugate}
f^p(M^*,\delta^*) = \sup_{M\in \m} 
\pro{M^*}{M} + \delta^* \det M -f(M),\quad \mbox{with}\ \ M^*\in\m \mbox{ and } \delta^*\in \re,
\end{equation}
and $f^{pp}(M) =(f^{p})^*(M,\det M),$ where   \[(f^{p})^*(M,\lambda) = \sup_{\substack{\delta^*\in \re \\ M^*\in \m}}\left(\pro{M}{M^*}+\delta^*\lambda -f^p(M^*,\delta^*)\right), \quad \mbox{for}\ \ M\in\m \mbox{ and } \lambda\in \re.
\]

\notelm{Modifique el parrafo anterior para agregar las correciones sugeridas. Debo volver a checarlo}In the proof of Theorems \ref{T3} and \ref{T2}, we strongly use the existence of the matrix $U_0$ 
given in Lemma \ref{origin}. By the  self-similarity of the incompatible cone and 
Lemma~\ref{planeQV}, we set $U_0$ as the center of coordinates.  The following lemma allow us to do that.     
\begin{proposition}\label{trans}
Let $f:\m \rightarrow \re$ and $g(M) = f(M-M_0)$, then $g^{pp}(M) = f^{pp}(M-M_0)$. 
\end{proposition}
\begin{proof}
From relation \eqref{fconjugate}, 
\[
\begin{array}{rl}
g^p(M^*,\delta^*) =& \sup_{M\in \m}\left(\pro{M^*}{M}+\delta^* \det M - f(M-M_0)\right)\\ 
=& \sup_{\bar{M}\in \m}\left(\pro{M^*}{\bar{M}+M_0}+\delta^* \det (\bar{M}+M_0) - f(\bar{M})\right).
\end{array} 
\]
The last expression for $g^p(M^*,\delta^*)$ is easily simplified by  Remark~\ref{det:prop}~\ref{det:prop3}. Hence 
\[
\begin{array}{rl}
g^p(M^*,\delta^*) =\pro{M^*}{M_0} + \delta^* \det M_0 + f^p(M^* + \delta^* \mathcal{S}M_0, \delta^*).
\end{array} 
\]
Moreover,  
\[
\begin{array}{rl}
(g^p)^*(M,\delta) 
=& \sup_{\substack{M^*\in \m \\ \delta^* \in \re}}\left\{\pro{M-M_0}{M^*} + \delta^*(\delta-\det M_0)- f^p(M^* + \delta^* \mathcal{S}M_0, \delta^*)\right\},\\ 
=& \sup_{\substack{\bar{M}\in \m \\ \delta^* \in \re}} \left\{ \pro{M-M_0}{\bar{M}} + \delta^*(\delta-\det M_0-\pro{M-M_0}{\mathcal{S}M_0}) - f^p(\bar{M}, \delta^*)\right\}.
\end{array}
\]
Thus, $(g^p)^*(M,\delta)=  (f^p)^*(M-M_0, \delta-\det M_0-\pro{M-M_0}{\mathcal{S}M_0})$. 
Finally, the affirmation follows from $g^{pp}(M)=(g^p)^*(M,\det M)$, and 
$\pro{M-M_0}{\mathcal{S}M_0} = -2\det M_0+\pro{\mathcal{S}M_0}{M}$.
\end{proof}

The next proposition determines the polyconvex conjugate function of $f:\m\rightarrow \re$ when it depends on the symmetric part of the argument.

\begin{proposition}\label{splitsup}
Let $f:\m \rightarrow \re$ be such that $f(M) = f(e(M))$ for every $M\in\msym$. Then, the polyconvex conjugate function of $f$ is given by 
\begin{equation}\label{f*sym}
f^p(M^*,\delta^*) = \begin{cases}
 \infty, & \mbox{if } \delta^*>0 \ \mbox{or } \delta^* = 0,\ w_{M^*}\neq 0,\\
 -\frac{w_{M^*}^2}{\delta^*} + \sup_{U\in \msym}\mathcal{L}(e(M^*),\delta^*, U), & \mbox{if } \delta^*<0, \\
 \sup_{U\in \msym}\mathcal{L}(e(M^*),0, U), & \mbox{if } \delta^* = w_{M^*} = 0, 
 \end{cases}
 \end{equation}
where $w_M\in \re$ satisfies $w(M)= w_M R$ and 

\begin{equation}\label{eqL}
\mathcal{L}(V,\delta, U):= \pro{V}{U} + \delta \det U-f(U).
\end{equation}

\end{proposition}
\begin{proof}
 From \eqref{fconjugate} and Remark \ref{det:prop}~\ref{det:prop2}, we get
\[
\begin{array}{rl}
f^p(M^*,\delta^*) =& \sup_{M\in \m} \left[\pro{M^*}{M} + \delta^* \det M - f\left(e(M)\right)\right]\\
=& \sup_{\substack{M\in \m}}\left[ \begin{array}{l} 2w_{M}w_{M^*} + \delta^* w_M^2  
+ \pro{e(M^*)}{e(M)} + \delta^* \det e(M) - f\left(e(M)\right)\end{array}\right].
\end{array}
\]
Since $w_M$ and $e(M)$ are independent, we conclude that
\begin{equation}\label{f*}
f^p (M^*, \delta^*)=\sup_{w\in \re}\left(2ww_{M^*}+ \delta^*w^2\right) + \sup_{U\in \msym}\mathcal{L}(e(M^*),\delta^*, U).
\end{equation}
It is not difficult to see that $f^p (M^*, \delta^*) = \infty$ if either $\delta^*>0$ or $\delta^* = 0,$ and $w_{M^*}\neq 0$. 
 Moreover, if $\delta^*<0$, the supremum in the first term on the right-hand side of \eqref{f*} is attained at $w=-w_{M^*}/\delta^*$. 
Thus 
\[
f^p (M^*, \delta^*)=-\frac{w_{M^*}^2}{\delta^*} + \sup_{U\in \msym}\mathcal{L}(e(M^*),\delta^*, U).
\]
Finally, the third condition in \eqref{f*sym} follows directly from \eqref{f*}, and this finishes the proof. 
\end{proof}

The following proposition characterizes the polyconvex biconjugate function of symmetric function, that is of the form $f(M)=f(e(M))$.

\begin{proposition}\label{fPPsym}
Let $f:\m \rightarrow \re$.  If $f(M) = f(e(M))$ for every $M\in \m$, then  $f^{pp}(M)=f^{pp}(e(M))$ and 
\[
f^{pp}(M) = \sup_{\substack{ V\in \msym, \delta^*\leq 0}} 
\left\{\pro{e(M)}{V} +  \delta^*\det e(M) - \sup_{U\in \msym} \mathcal{L}(V,\delta^*,U) \right\}.
\]
\end{proposition}

\begin{proof}
Since $f(M) = f(e(M))$,  $f^p(M^*,\delta^*)$ is given by \eqref{f*sym}. Hence, 
\[
(f^{p})^*(M,\delta)=\max\{g_1(M),g_2(M,\delta)\}
\]
where 
\begin{eqnarray*}
&g_1(M)=\sup_{ V\in \msym}\left\{ \pro{e(M)}{V}  - \sup_{U\in \msym} \mathcal{L}(V,0,U) \right\}, \text{ and } \\
&g_2(M,\delta)=\sup_{\substack{w\in\re, \delta*< 0\\ V\in \msym}}\left\{
 \pro{e(M)}{V} +2w_M w  + \frac{w^2}{\delta^*} + \delta \delta^* - \sup_{U\in \msym} \mathcal{L}(V,\delta^*,U)\right\},
\end{eqnarray*}
with $\mathcal{L}(V,\delta^*,U)$ as in \eqref{eqL}.  Now, the function
\[
w\mapsto  \pro{e(M)}{V} + 2w_M w  + \frac{w^2}{\delta^*} + \delta \delta^* - \sup_{U\in \msym} \mathcal{L}(V,\delta^*,U)
\] 
is concave in $w$ for every $\delta^*<0$. Thus, 
\[
\begin{array}{rl}
    2w_M w  + \frac{w^2}{\delta^*} + \pro{e(M)}{V} + \delta \delta^* - \sup_{U\in \msym} \mathcal{L}(V,\delta^*,U) \leq & \\ & \hspace{-2.5cm} \pro{e(M)}{V} + (\delta-w_M^2) \delta^* - \sup_{U\in \msym} \mathcal{L}(V,\delta^*,U),
\end{array}
\]
with equality at $w=-\delta^*w_M$. Therefore,
\[
(f^{p})^*(M,\delta) = \sup_{\substack{ V\in \msym,  \delta^*\leq 0}}
\left\{ \pro{e(M)}{V} + (\delta-w_M^2) \delta^* - \sup_{U\in \msym} \mathcal{L}(V,\delta^*,U)\right\}.
\]
Next, we notice that $w^2_M = \det w(M)$.  Thus, since  $f^{pp}(M)= (f^{p})^*(e(M),\det e(M))$, 
by \cref{det:prop}~\ref{det:prop2} we obtain
\[
\begin{array}{rl}
f^{pp}(M) = & \sup_{\substack{ V\in \msym, \delta^*\leq 0}} 
\left\{\pro{e(M)}{V} +  \delta^*\det e(M) - \sup_{U\in \msym} \mathcal{L}(V,\delta^*,U) \right\}.
\end{array}
\]
Hence, $f^{pp}(M)=f^{pp}(e(M))$ as claimed. 
\end{proof}

 \section{The polyconvex conjugate and biconjugate functions of \texorpdfstring{$f_C$}{}} \label{sec:7}

In this section we specialize our results to a particular function $f_C$ that we use in our proofs. Let $f:\m\rightarrow \re$ be such that 
\begin{equation}\label{eq:fc}
f_C(M) = \chi_{\bar{B}}(e(M))\left(|\det e(M)| + \left|\pro{C}{M}\right|\right),
\end{equation}
where  $C\in \msym$ with negative determinant, and 
\[
\chi_{\bar{B}}(U) =
\begin{cases}
1 & \text{ if } U\in \msym \setminus L^e(\mathcal{U})\\
0 & \text{otherwise.}
\end{cases}
\]

The following lemmas will be used to determine the polyconvex conjugate and  biconjugate functions of $f_C$ given as in \cref{eq:fc}. Notice that $f_C$, restricted to $\Aff(\mathcal{U})$, vanishes on $L^e(\mathcal{U})$. 
\todoin{From now on, we denote by $\mathcal{V}$ any  three-well set where $\Aff(\mathcal{V})$ is a two-dimensional subspace, and by $\mathcal{U}$ any three-well set where $\Aff(\mathcal{U})$ has codimension one. 
}
 
 \begin{lemma}\label{supOut}
Let $\mathcal{V} = \{V_1,V_2,V_3\}\subset \msym$  such that $\det(V_2)=\det(V_3)=0$, and $\Aff(\mathcal{V})\subset\msym$ is a two-dimensional subspace. Also, assume that either $\mathcal{V}$ is type one, or $\mathcal{V}$ is type two and $0\in C(\mathcal{V})$. 
Additionally, let $f$ be defined as in \cref{eq:fc} for a fix $C\in \msym$ with negative determinant, and 
$\mathcal{L}(V,\delta, U)$  be defined as in \cref{eqL}. Then



 
\begin{equation}\label{sup11}
\sup_{U\in \msym}\mathcal{L}(V,\delta^*, U) = \begin{cases}
\infty & \mbox{if } (V,\delta^*)\in\mathcal{N}^c,\\
\max_{U\in \mathcal{V}\cup\{0\}} (k\pro{C}{U}+\delta^*\det U) & \mbox{if } \ 
(V,\delta^*)\in\mathcal{N},
\end{cases}
\end{equation}
\todoin{where $\mathcal{N}$ and  $\mathcal{N}^c$ are subsets of $\msym\times (-\infty,0]$  given by
\[
\mathcal{N} = \left\{ (kC,\delta^*) \  \middle| \ -1<k<1, \ -1\leq\delta^*\leq 0\right\} \ \ \mbox{and} \ \ \mathcal{N}^c = \left(\msym\times (-\infty,0]\right) \setminus \mathcal{N} .
\]
}

\end{lemma}

\begin{proof} We observe that $\Aff(\mathcal{V})=\Pi_Q(0)$, for some $Q\in \msym$ with negative determinant since  there are a compatible and an incompatible pair of wells in $\Aff(\mathcal{V})$ (see \cref{rem:det<}).  We divide the proof into two parts corresponding to the following claims. First, we claim that 
\begin{equation}\label{SymSup}
\sup_{U\in \msym \setminus L^e(\mathcal{V})}\mathcal{L}(V,\delta^*, U) = \begin{cases}
0 & \mbox{if } (V,\delta^*)\in \mathcal{N},\\ 
\infty & \mbox{otherwise.} 
\end{cases}
\end{equation}
Second, we claim that if $(V,\delta^*)\in \mathcal{N}$, then
\begin{equation}\label{SymSup2}
\sup_{U\in L^e(\mathcal{V})}\mathcal{L}(kC,\delta^*, U) = \begin{cases}
\max_{U\in \mathcal{V}} (k\pro{C}{U}+\delta^*\det U) & \mbox{if }  \ \mathcal{U} \ \ \mbox{is type one},\\
\max_{U\in \mathcal{V}\cup\{0\}} (k\pro{C}{U}+\delta^*\det U) & \mbox{if } \ \mathcal{U} \ \ \mbox{is type two.}\end{cases}
\end{equation}

Assuming the claims, we finish the proof. 
Let 
$L:=\sup\{\mathcal{L}(V,\delta^*,U)\,|\, U\in \msym \}$ for each $(V,\delta^*)\in \msym\times(-\infty,0]=\mathcal{N}\cup\mathcal{N}^c$. Hence, if $(V,\delta^*)\in \mathcal{N}^c$, then $L=\infty$ by \cref{SymSup}, and if $(V,\delta^*)\in \mathcal{N}$, then 
\[
L=\max\left\{ \sup\{\mathcal{L}(V,\delta^*,U)\,|\, U\in L^e(\mathcal{V}) \}, \sup\{\mathcal{L}(V,\delta^*,U)\,|\, U\in \msym\setminus L^e(\mathcal{V}) \}\right\}.
\]
Therefore, the second part of \cref{sup11} is a consequence of \cref{SymSup} and \cref{SymSup2}, and the results follows. 

Next, we divide the proof of the claims in two steps:\\
{\bf Step 1:} We will prove \cref{SymSup}. Since $U\notin L^e(\mathcal{V})$, then   
\[\mathcal{L}(e(M^*),\delta^*,U) = \min\{(\delta^*-1) \det U, (\delta^*+1) \det U \} + \min\{\pro{e(M^*)-C}{U}, \pro{e(M^*)+C}{U}\}.\]
We notice that 
\[\mathcal{N}^c=\{(V,\delta^*)\in \msym\times(-\infty,0]\ \  |\ \ \delta^*<-1\ \mbox{ or }\ V \nparallel C \ \mbox{ or }\ V=kC,\ |k|>1\}.
\]
Now, we inspect each option for $(V,\delta^*)\in \mathcal{N}^c$.
\begin{enumerate}
\item[(i)] Assume that $\delta^*<-1$. 
Since $(\delta^*-1),(\delta^*+1) < 0$, we choose 
$U(t)=t\bar{U}$ for a fix $\bar{U}\in\msym \setminus L^e(\mathcal{V})$ with $\det \bar{U}<0$, and  $t\in\mathbb{R}$. We find
\[
\mathcal{L}(e(M^*),\delta^*, U(t)) = t^2(\delta^*-1)\det \bar{U} + t\min\{\pro{e(M^*)-C}{\bar{U}}, \pro{e(M^*)+C}{\bar{U}}\}.
\]
By letting $t\rightarrow \infty$ the claim follows.
\medskip 

\item[(ii)] Now, we assume $e(M^*)$ and $C$ linearly independent. Thus, there exists  $V\in\msym$ with $\pro{V}{C}=0$ such that
$e(M^*) = \alpha V + kC$ for some $\alpha, k\in \re$. Thus,  $V\in \Pi_C(0)$. Due to $\det C<0$ and Lemma \ref{origin}, the boundary of  
$\Pi_C(0)\cap\mathcal{C}_{0}$ consist of two nonparallel rank-one lines, say \linebreak \hfill $\ell_1=\{tu_1\otimes u_1 | t\in \re\}$ and $\ell_2=\{tu_2\otimes u_2 | t\in \re\}$ for some $u_1,u_2\in S^1$. Therefore, we choose $\bar{U}\in\{u_1\otimes u_1,\, -u_1\otimes u_1,\, u_2\otimes u_2,\, -u_2\otimes u_2\}$  
such that $\alpha\pro{V}{\bar{U}}>0$. By construction, $\det \bar{U} = 0$ and $\pro{C}{\bar{U}}=0$, hence 
$\mathcal{L}(e(M^*),\delta^*, t\bar{U}) =  t\alpha\pro{V}{\bar{U}}$ and by letting $t\rightarrow \infty$ we get the result. 
\medskip 

\item[(iii)] If $e(M^*) =  kC$  with $|k|>1$, then we choose $\bar{U}\in \mathcal{C}_{0}\setminus \Pi_C(0)$ such that 
$k\pro{C}{\bar{U}}>0$  and  $\det \bar{U}=0$. Hence, 
$\mathcal{L}(e(M^*),\delta^*, t\bar{U}) = \min\{t(k-1)\pro{C}{\bar{U}}, t(k+1)\pro{C}{\bar{U}}\}$. Clearly, $(k-1)\pro{C}{\bar{U}}>0,$ and $(k+1)\pro{C}{\bar{U}}>0$, so by letting $t\to \infty$ the result follows.
\end{enumerate}
From the last three steps, we conclude that $\sup\{\mathcal{L}(V,\delta^*, U)\,|\, U\in \msym \setminus L^e(\mathcal{V})\}=\infty,$ if $(V,\delta^*)\in\mathcal{N}^c$. 

Now, we assume $(V,\delta^*)\in \mathcal{N}$, and we prove that $\sup\{\mathcal{L}(e(M^*),\delta^*,U)\,|\, U\in \msym\setminus L^e(\mathcal{V})\} = 0$. Indeed, we have that 
\begin{equation}\label{sup0}
\mathcal{L}(kC,\delta^*,U) = \min\{(\delta^*-1) \det U, (\delta^*+1) \det U \} + \min\{(k-1)\pro{C}{U}, (k+1)\pro{C}{U}\}\leq 0.
\end{equation}
The last inequality follows since both terms on the right hand side of \eqref{sup0} are non-positive. Thus if  $\bar{U}\in \mathcal{C}_{0}\setminus \Pi_{Q}(0)$ is fixed, then 
we get 
\[
\mathcal{L}(kC,\delta^*,t\bar{U}) = \min\{t(k-1)\pro{C}{\bar{U}}, t(k+1)\pro{C}{\bar{U}}\} \rightarrow 0, \quad \mbox{as} \quad t\rightarrow 0.
\]
Therefore, the affirmation follows and the claim (i.e. \cref{SymSup}) is proved.
\medskip

{\bf Step 2:} Now, we consider the optimization of ${\mathcal L}$ on $U\in L^e(\mathcal{V})$ (i.e. \cref{SymSup2}). By \cref{origin}, for every $U\in L^e(\mathcal{V})\subset \Pi_Q(0)$ there exist  $\xi,\, \eta \in \re$ such that $U = \pv{\xi}{\eta},$ and 
\begin{equation}\label{simpSup}
\mathcal{L}(kC,\delta^*,U)=k\pro{C}{U} +\delta^*\det U = \xi \eta\, \delta^*|n\times a|^2 +\xi Ca^\perp \cdot a^\perp+ \eta Cn^\perp \cdot n^\perp.
\end{equation}
The level sets of \eqref{simpSup} as function of $\xi$ and $\eta$ are hyperbolae (or straight lines). 
Hence, the supremum of $\mathcal{L}(kC,\delta^*,U)$ in $L^e(\mathcal{V})$ is attained at some point $\bar{U}$ on the relative boundary of $L^e(\mathcal{V})$, denoted by $\rpartial L^e(\mathcal{V})$. By Proposition \ref{P:lamconv},  $\rpartial L^e(\mathcal{V})$ depends on the type of $\mathcal{V}$, namely 
\[
\rpartial L^e(\mathcal{V}) = \begin{cases}
L^e(\mathcal{V}) & \mbox{if }\mathcal{V} \mbox{ is type one,}\\
 C(\{V_1,\,V_2\}) \cup C(\{V_2,\,0\}) \cup C( \{0,\,V_3\}) \cup C(\{V_3,\,V_1\}) & \mbox{if }\mathcal{V} \mbox{ is type two}.
 \end{cases}
\]
First, we assume that $\mathcal{V}$ is type two. If $\bar{U}\in C(\{V_2,\,0\}) \cup C( \{0,\,V_3\})$, then $\det\bar{U}=0$ and the equation \eqref{simpSup} becomes linear. Hence, the maximum of $\mathcal{L}(kC,\delta^*,U)$ in $C(\{V_2,\,0\}) \cup C( \{0,\,V_3\})$ is attained at either $0,$ $V_2$ or $V_3$. Now, if  $U\in C(\{V_1,\,V_2\})$, then 
$U = \pv{[\lambda\xi+(1-\lambda)\gamma]}{\lambda \eta}$ 
 for some $\lambda\in [0,1]$ due to \cref{cor:origin}. So, \cref{simpSup} becomes a polynomial of degree two in $\lambda$. It is readily seen that 
\[
\frac{d^2\mathcal{L}(kC,\delta^*, U(\lambda))}{d\lambda^2} = 2\delta^* \eta (\xi-\gamma)|n\times a|^2\geq 0.
\]
Thus, either $\mathcal{L}(kC,\delta^*, U(\lambda))$ is linear in $\lambda$ if $\xi=\gamma$ or it is quadratic in $\lambda$ with a minimum if $\xi\neq\gamma$. 
Therefore, the maximum value of $\mathcal{L}(kC,\delta^*, U)$ on $C(\{V_1,\,V_2\})$ 
is always attained at the extremal points  $\{V_1,V_2\}$. With a similar argument, the supremum of $\mathcal{L}(kC,\delta^*, U)$ on $C(\{V_1,V_3\})$ is attained on either $V_1$ or $V_3$. Therefore, we conclude that 
the supremum of $\mathcal{L}(kC,\delta^*,U)$ in $L^e(\mathcal{V})$  is attained at some matrix $\bar{U}\in\{0,V_1,V_2,V_3\}$ as claimed. 
\medskip 

Second, we assume that $\mathcal{V}$ is of type one. In this case,  the maximum value of $\mathcal{L}(kC,\delta^*,U)$ is attained either in $V_1$ or in $C(\{V_2,V_3\})$. If the second option happens,  $d^2\mathcal{L}/d\lambda^2=-2\delta^*\gamma\zeta|n\times a|^2\geq 0$, and by an   analogous argument, the result follows straight forward. 

Therefore, the supremum of $\mathcal{L}(kC,\delta^*,U)$ in $L^e(\mathcal{V})$  is attained at some matrix $\bar{U}\in\{V_1,V_2,V_3\}$ and \cref{SymSup2} follows. 
Now the proof is completed. 
\end{proof}

The last result of this section contains explicit bounds for $\Ker f^{pp}\subset \m$, the zero-level set of the polyconvex biconjugate function of $f:\m\rightarrow \re$. 
From now on, if $g:\m\rightarrow \re$ is any function of the form $g(M)=g(e(M))$ for every $M\in \msym$, we denote the zero-level set of $g|_{\msym}$ as
\[
\kers g:=\{V\in \msym\,|\, g(V)=0\}.
\] 
In this type of functions it readily follows that $\Ker g = \{M\in \m\, |\, e(M)\in \kers g\}$.
In particular, if we study $\Ker f^{pp}$, we can focus on $\kers f^{pp}$ and no information is lost.

\begin{lemma}\label{eqboundary} 
Let $\mathcal{V} = \{V_1,V_2,V_3\}\subset \msym$ be either type one or two 
such that $\Aff(\mathcal{V})\subset \msym$ is a two-dimensional subspace with $\det V_2=\det V_3=0$. Let $C\in \msym$ with negative determinant and define $f(M)$ as in \cref{eq:fc}. We have that:
\begin{enumerate}
\item If $C$ is normal to $\Aff(\mathcal{V})$, then $f^{pp}(M)\geq 0$ for every $M\in \m$, and 
\begin{equation}\label{kerCP}
\Ker f^{pp} = \left\{M\in \m\,\middle | \, e(M)\in \Aff(\mathcal{V}),\, \det V_1 \leq \det e(M) \right\}.
\end{equation}
\item If $0\in C(\mathcal{V})$, $C$ satisfy $\pro{C}{V_1}\leq 0< \pro{C}{V_3}=\pro{C}{V_2}$, and 
\begin{enumerate}[label=(\roman*)]
\item\label{lem9:itemb2} $\mathcal{V}$ is type one, then
\begin{equation*}
\begin{split}
\Ker f^{pp}\subset \{ M \in \m| \pro{C}{e(M)-V_2} \leq 0  \leq &\pro{C}{e(M)-V_1}, \\
&0 \leq \det e(M) \ \text{ and } \ 0\leq \hbar(M)\},  
\end{split}
\end{equation*}

\item\label{lem9:itemb1} $\mathcal{V}$ is type two, then
\begin{equation*}
\begin{split}
\Ker f^{pp}\subset \{ M \in \m| \pro{C}{e(M)-V_2} \leq 0  \leq &\pro{C}{e(M)-V_1}, \\
&\det V_1 \leq \det e(M) \ \text{ and } \ 0\leq \hbar(M)\}.  
\end{split}
\end{equation*}
\end{enumerate}

where 
\begin{equation*}
\hbar(M) = \pro{C}{e(M)-V_2}\det V_1 - \pro{C}{V_1-V_2}\det e(M).
\end{equation*}
Moreover, $f^{pp}$ is a non negative real-valued function and  $\mathcal{V} \subset \kers f^{pp}$.
\end{enumerate}
\end{lemma}

\begin{proof} 
From \cref{fPPsym} and \cref{supOut} we have that 
\begin{equation}\label{fPPspec}
f^{pp}(M) =  \sup_{(k,\delta^*)\in\mathcal{E}} \left( k\pro{e(M)}{C}+\delta^*\det e(M) 
-\max_{U\in \mathcal{V}\cup\{0\}} (k\pro{C}{U}+\delta^*\det U)  \right),
\end{equation}
where $\mathcal{E} := [-1,1]\times [-1,0]$. 
Our goal is to estimate outer bounds on $\kers f^{pp}$.

First, we assume that $C$ is normal to $\Aff(\mathcal{V})$.  Under these assumptions and since $0\in\Aff(\mathcal{V})$,  \cref{fPPspec} becomes
\[
f^{pp}(V)=\sup_{(k,\delta^*)\in\mathcal{E}} (k\pro{V}{C}+\delta^*\left(\det V -\det V_1\right)).
\]
By computing the supremum, we get
\begin{equation}\label{fPPspecQ}
f^{pp}(V) = |\pro{V}{C}|+ \max\{0,\det V_1-\det V\}.
\end{equation} 
The first term on the right hand side of \eqref{fPPspecQ} penalizes the distance to the plane $\Aff(\mathcal{V})$, 
meanwhile the second term is an in-plane condition. Moreover, $f^{pp}$ is non negative and  
\[
\kers f^{pp}= \left\{V\in \msym\,\middle | \, V\in \Aff(\mathcal{V}),\, \det V_1 \leq \det V \right\},
\]
as claimed in \cref{kerCP}.  
\medskip 

Second, we assume $\pro{C}{V_1}\leq 0<\pro{C}{V_2}=\pro{C}{V_3}$. A simple computation lead us to
\[
k\pro{C}{V_2}> \delta^* \det V_1 + k\pro{C}{V_1}, \quad \mbox{for } \quad \frac{\det V_1}{\pro{C}{V_2-V_1}}\delta^* \leq k \leq 1.
\]
This observation implies that 
\begin{equation}\label{eq:maxL}
\max_{i=0,\, 1,\, 2,\, 3}\{\delta^* \det V_i + k\pro{C}{V_i} \}= \begin{cases}
 k\pro{C}{V_2}, & \mbox{if } \ \ \frac{\det V_1}{\pro{C}{V_2-V_1}}\delta^* \leq k \leq 1,\\
 \delta^*\det V_1 + k\pro{C}{V_1}, & \mbox{otherwise.}
 \end{cases}
\end{equation}
Now, we define the sets $\mathcal{D} = \left\{(k,\delta^*)\in\mathcal{E}\,:\,  \delta^*\det V_1/\pro{C}{V_2-V_1} \leq k \leq 1\right\},$ 
and its complement 
\[
\mathcal{E}\setminus \mathcal{D} =\left\{(k,\delta^*)\in \mathcal{E}\,\middle|\, -1\leq k\leq \min \left\{1,\frac{-\delta^*\det V_1}{\pro{C}{V_1-V_2}}\right\}, \ \ -1\leq \delta^*\leq 0\right\},
\]
relative to $\mathcal{E}$, see Figure \ref{fig:div1}. 
Let $\chi_{\mathcal{D}}:\mathcal{E}\rightarrow \{0,1\}$ and 
$\chi_{\mathcal{E}\setminus \mathcal{D}}:\mathcal{E}\rightarrow \{0,1\}$ be the characteristic functions of 
$\mathcal{D}$ and $\mathcal{E}\setminus \mathcal{D}$, respectively.
From  \cref{fPPspec} and \cref{eq:maxL}, we get
\[
f^{pp}(V) = \sup\{\bar{f}_1(V,k,\delta^*) + \bar{f}_2(V,k,\delta^*)\, | \, (k,\delta^*)\in \mathcal{E}\},
\]
where 
\begin{equation*}
\begin{cases} 
\bar{f}_1(V,k,\delta^*) =  \chi_{\mathcal{D}}\left ( k\pro{C}{V-V_2}+\delta^*\det V \right), \\ 
\bar{f}_2(V,k,\delta^*) =\chi_{\mathcal{E}\setminus \mathcal{D}} \left( k\pro{C}{V-V_1}+\delta^*(\det V-\det V_1)\right).
\end{cases} 
\end{equation*}
\begin{figure}[h]
\centering
\begin{tikzpicture}[thick,scale=0.65, every node/.style={transform shape}]
\begin{scope}
\filldraw[color =white, fill=black!60] (0,0) -- (3.5,-2.75) -- (3.5,0) -- cycle;
\filldraw[color=white, fill=black!30] (0,0) -- (-3.5,0) -- (-3.5,-3.5) --(3.5,-3.5) -- (3.5,-2.75) -- cycle;
\draw (0,0) -- (3.5,-2.75);
\draw[thick, ->] (-4,0) -- (4,0);
\node at (4,-.4) {$k$};
\draw[thick, ->] (0,-4) -- (0,.5);
\node at (0.4,.5) {$\delta^*$};
\draw[dotted,thick] (-3.5,0) -- (-3.5,-3.5) -- (3.5,-3.5) -- (3.5,0); 
\draw[thick] (-0.2,-3.5) -- (0.2,-3.5) node [anchor=north]{-1};
\draw[thick] (-3.5,-0.2) -- (-3.5,0.2) node [anchor=south]{-1};
\draw[thick] (3.5,-0.2) -- (3.5,0.2) node [anchor=south]{1};
\node at (2.5,-1) {\large $\mathcal{D}$};
\node at (-1.1,-2.5) {\large $\mathcal{E}\setminus \mathcal{D}$};

\end{scope}

\begin{scope}[shift={(9,0)}]
\filldraw[color =white, fill=black!60] (0,0) -- (2.75, -3.5) -- (3.5,-3.5) -- (3.5,0) -- cycle;
\filldraw[color=white, fill=black!30] (0,0) -- (-3.5,0) -- (-3.5,-3.5) --(2.75,-3.5) -- cycle;
\draw (0,0) -- (2.75,-3.5);
\draw[thick, ->] (-4,0) -- (4,0);
\node at (4,-.4) {$k$};
\draw[thick, ->] (0,-4) -- (0,.5);
\node at (0.4,.5) {$\delta^*$};
\draw[dotted,thick] (-3.5,0) -- (-3.5,-3.5) -- (3.5,-3.5) -- (3.5,0); 
\draw[thick] (-0.2,-3.5) -- (0.2,-3.5) node [anchor=north]{-1};
\draw[thick] (-3.5,-0.2) -- (-3.5,0.2) node [anchor=south]{-1};
\draw[thick] (3.5,-0.2) -- (3.5,0.2) node [anchor=south]{1};
\node at (2.5,-1) {\large $\mathcal{D}$};
\node at (-1.1,-2.5) {\large $\mathcal{E}\setminus \mathcal{D}$};

\end{scope}

\end{tikzpicture}
\caption{Two possible subdivisions of the set $\mathcal{E}$ when $\mathcal{V}$ is of type two. Left and right images correspond to  the cases  $\pro{C}{V_2-V_1}>-\det V_1$ and $\pro{C}{V_2-V_1}<-\det V_1$, respectively. 
}
\label{fig:div1}
\end{figure}
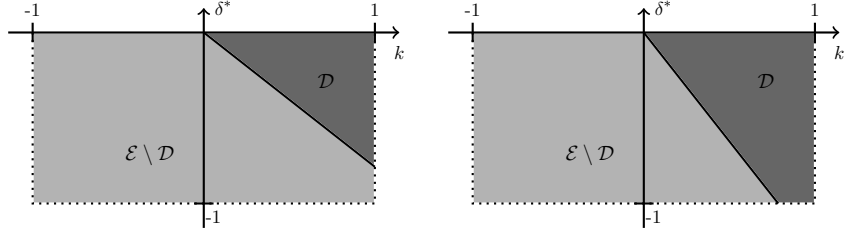
Now, $\mathcal{D}$ and $\mathcal{E}\setminus \mathcal{D}$ are disjoint sets, hence 
\[
f^{pp}(V) = \max\{f_1(V),f_2(V)\}, \quad  \mbox{where} \quad \begin{cases}
f_1(V) = \sup_{(k,\delta^*)\in \mathcal{D}}\bar{f}_1(V,k,\delta^*), \\ f_2(V) =\sup_{(k,\delta^*)\in \mathcal{E}\setminus \mathcal{D}}\bar{f}_2(V,k,\delta^*). 
 \end{cases}
\]
Notice that $\bar{f}_1$ is linear in $(k,\delta^*)$. Thus, its supremum in $\mathcal{D}$ is attained on $\partial \mathcal{D}$ for each $V\in \msym$. Analogously, the supremum of $\bar{f}_2$ in 
$\mathcal{E}\setminus \mathcal{D}$ is attained on $\partial (\mathcal{E}\setminus \mathcal{D})$.
The common boundary between $\mathcal{D}$ and $\mathcal{E}\setminus \mathcal{D}$ depends on  $\det V_1$, particularly on its sign (see Figure \ref{fig:div1}). Therefore, the computation of $f_1$ and $f_2$ also depends on $V_1$.  Moreover, since $0\in C(\mathcal{V})$ and $\det V_2=\det V_3=0$, \cref{cor:origin} implies that $\det V_1$ is positive or negative if $\mathcal{V}$ is of type one or two, respectively.  

Now, we assume $\mathcal{V}$ is type two. Thus, $\det V_1\le 0$ and we divide the computation of $f_1$ into four cases:
\begin{enumerate}[label=(\alph*)]
\item Assume  $0<\pro{C}{V-V_2}$ and $0<\det V$, then the supremum is attained at $(k,\delta^*) = (1,0),$ 
and $f_1(V) = \pro{C}{V-V_2}> 0$. 
\item Assume $0<\pro{C}{V-V_2}$ and $\det V\leq 0$. If $\det V_1 <0$, then $\bar{f}_1$'s supremum is attained at $(k,\delta^*) = \left(1,\max\left\{-1,-\pro{C}{V_1-V_2}/\det V_1 \right\}\right)$, and if 
$\det V_1 = 0$ then $\bar{f}_1$'s supremum  is attained at $(k,\delta^*)=(1,-1)$. In either case,
\[
f_1(V) = \pro{C}{V-V_2}-\frac{1}{\max\left\{1,\frac{\det V_1}{\pro{C}{V_1-V_2}} \right\}}\det V>0.
\]
\item Assume $\pro{C}{V-V_2}\leq 0$ and $\det V\leq 0$. In this case, the 
supremum of $\bar{f}_1$ is attained on the line segment 
\[
\left\{\left(-\frac{\det V_1}{\pro{C}{V_1-V_2}}\delta^*,\delta^*\right)\in \mathcal{E}\,\middle | \, \delta^*\in [-1,0] \right\},
\] 
and the function $\bar{f}_1$ evaluated at this segment is  
\[
\bar{f}_1\left(V, -\frac{\det V_1}{\pro{C}{V_1-V_2}}\delta^*,\delta^*\right) = -\frac{\hbar(V)}{\pro{C}{V_1-V_2}}\delta^*,
\]
where $\hbar(V) = \pro{C}{V-V_2}\det V_1 - \pro{C}{V_1-V_2}\det V$.
The last equation is linear in $\delta^*$ and  
increasing or decreasing depending on the sign of 
$\hbar(V)/\pro{C}{V_1-V_2}$. 
Hence, if $\hbar(V)$ is non negative, $f_1(V) = 0$, meanwhile for $\hbar(V)<0$, 
\[
 f_1(V) = \min\left\{1,\frac{\det V_1}{\pro{C}{V_1-V_2}}\right\}\pro{C}{V-V_2} - \frac{1}{\max\left\{1,\frac{\det V_1}{\pro{C}{V_1-V_2}}\right\}}\det V. 
\]
Inserting the definition of $\hbar(V)$ in the last equation we get
\[
f_1(V) =\frac{\hbar(V)}{p}>0,\quad \mbox{where} \quad p=\begin{cases}\det V_1, &  \mbox{if} \quad 1 \leq  \det V_1/\pro{C}{V_1-V_2},\\
\pro{C}{V_1-V_2}, & \mbox{otherwise}.  
\end{cases}
\]
\item Let $\pro{C}{V-V_2}\leq 0$ and $0< \det V$. Then, the supremum of $f_1$ is attained at  $(k,\delta^*) = (0,0)$ and $f_1(V) = 0$. 
In this case, we also have that  $\hbar(V)>0$. 
\end{enumerate}

Summarizing the above cases, we conclude that
\[
\kers f_1 =\{ V \in \msym \,| \,  \pro{C}{V-V_2}\leq 0  \mbox{ and either } \det V \leq 0\leq \hbar(V),\, \mbox{ or } 0 < \det V\},
\]
and we easily conclude 
\begin{equation}\label{ker1}
\kers f_1 \subset\{ V \in \msym \,| \,  \pro{C}{V-V_2}\leq 0 \leq \hbar(V)\}.
\end{equation}

Next,  we compute $f_2$. Again, we also have four possible further cases:
 \begin{enumerate}[label=(\alph*)]
\item If $\pro{C}{V-V_1}< 0$ and $0\leq \det V-\det V_1$, then the supremum  of $f_2$ is attained at  $(k,\delta^*) = (-1,0)$
and $f_2(V) = -\pro{C}{V-V_1} > 0$.
\item If $\pro{C}{V-V_1} \leq 0$ and $\det V-\det V_1<0$,  then the supremum is attained at  $(k,\delta^*) = (-1,-1)$
 and $f_2(V) = -\pro{C}{V-V_1} -\det V+\det V_1 > 0.$
\item If $0<\pro{C}{V-V_1}$ and $\det V-\det V_1< 0$,  then the supremum  of $f_2$ is attained at 
$(k,\delta^*) = \left(\min\left\{1,\det V_1/\pro{C}{V_1-V_2} \right\},-1\right)$ and 
\[
f_2(V) = \min\left\{1,\frac{\det V_1}{\pro{C}{V_1-V_2}} \right\}\pro{C}{V-V_1}-\det V+\det V_1>0.
\]
\item If $0\leq \pro{C}{V-V_1}$ and $0\leq \det V-\det V_1$,  then the supremum is attained at a point on the segment
\[
\left\{\left(\frac{-\det V_1}{\pro{C}{V_1-V_2}}\delta^*,\delta^*\right)\in \mathcal{E}\,\middle | \, \delta^*\in [-1,0] \right\},
\]
and 
\[
\bar{f}_2\left(V, \frac{-\det V_1}{\pro{C}{V_1-V_2}}\delta^*,\delta^*\right) = \frac{-\hbar(V)}{\pro{C}{V_1-V_2}}\delta^*.
\]
By linearity in $k$, if $-\hbar(V)/ \pro{C}{V_1-V_2}\geq 0$, then $f_2(V) = 0$. Meanwhile, for $\hbar(V)<0$,
 \[
f_2(V) = \min\left\{1,\frac{\det V_1}{\pro{C}{V_1-V_2}}\right\}\pro{C}{V-V_1} 
- \frac{1}{\max\left\{1,\frac{\det V_1}{\pro{C}{V_1-V_2}}\right\}}(\det V-\det V_1). 
\]
From the definition of $\hbar(V)$, we get
\[
f_2(M) =\frac{\hbar(V)}{p}>0,\quad \mbox{where} \quad p=\begin{cases}\det V_1, &  \mbox{if} \quad 1 \leq  \det V_1/\pro{C}{V_1-V_2},\\
\pro{C}{V_1-V_2}, & \mbox{otherwise}.  \end{cases}
\]
\end{enumerate}

Hence, from the above four cases, we conclude that 
\begin{equation}\label{ker2}
\kers f_2=\{ V \in \msym \,| \, 0 \leq \pro{C}{V-V_1}, \, \det V_1 \leq \det V, \,  0\leq \hbar(V)\}.
\end{equation}

Now, from the above analysis, we conclude that $f^{pp}_C$ is nonnegative and  $\kers\, f^{pp} = \kers\, f_1 \cap \kers\, f_2$.
Hence, by \eqref{ker1} and \eqref{ker2}, we get 
\begin{equation}\label{kertype2}
\kers f^{pp}\subset\{ V \in \msym \,| \, \pro{C}{V-V_2} \leq 0 \leq \pro{C}{V-V_1}, \, \det V_1 \leq \det V, \,  0\leq \hbar(V)\}, 
\end{equation}
that is \cref{lem9:itemb1}. 
A similar relation is obtained when 
%
$\mathcal{V}$ is type one. 
Under this assumption, we can assume $\det V_1>0$ and  the functions $f_1$ and $f_2$ can be computed analogously.
 The main observation is that  $\mathcal{D}$ and $\mathcal{E}\setminus \mathcal{D}$ are different due to the sign of $\det V_1$, see \cref{fig:div2}. 
 \begin{figure}[h]
\centering
\begin{tikzpicture}[thick,scale=0.65, every node/.style={transform shape}]
\begin{scope}
\filldraw[color =white, fill=black!30] (0,0) -- (-3.5,-2.75) -- (-3.5,0) -- cycle;
\filldraw[color=white, fill=black!60] (0,0) -- (3.5,0) -- (3.5,-3.5) --(-3.5,-3.5) -- (-3.5,-2.75) -- cycle;
\draw (0,0) -- (-3.5,-2.75);
\draw[thick, ->] (-4,0) -- (4,0);
\node at (4,-.4) {$k$};
\draw[thick, ->] (0,-4) -- (0,.5);
\node at (0.4,.5) {$\delta^*$};
\draw[dotted,thick] (-3.5,0) -- (-3.5,-3.5) -- (3.5,-3.5) -- (3.5,0); 
\draw[thick] (-0.2,-3.5) -- (0.2,-3.5) node [anchor=north]{-1};
\draw[thick] (-3.5,-0.2) -- (-3.5,0.2) node [anchor=south]{-1};
\draw[thick] (3.5,-0.2) -- (3.5,0.2) node [anchor=south]{1};
\node at (-2.5,-1) {\large $\mathcal{E}\setminus \mathcal{D}$};
\node at (1.1,-2.5) {\large $\mathcal{D}$};

\end{scope}

\begin{scope}[shift={(9,0)}]
\filldraw[color =white, fill=black!30] (0,0) -- (-2.75, -3.5) -- (-3.5,-3.5) -- (-3.5,0) -- cycle;
\filldraw[color=white, fill=black!60] (0,0) -- (3.5,0) -- (3.5,-3.5) --(-2.75,-3.5) -- cycle;
\draw (0,0) -- (-2.75,-3.5);
\draw[thick, ->] (-4,0) -- (4,0);
\node at (4,-.4) {$k$};
\draw[thick, ->] (0,-4) -- (0,.5);
\node at (0.4,.5) {$\delta^*$};
\draw[dotted,thick] (-3.5,0) -- (-3.5,-3.5) -- (3.5,-3.5) -- (3.5,0); 
\draw[thick] (-0.2,-3.5) -- (0.2,-3.5) node [anchor=north]{-1};
\draw[thick] (-3.5,-0.2) -- (-3.5,0.2) node [anchor=south]{-1};
\draw[thick] (3.5,-0.2) -- (3.5,0.2) node [anchor=south]{1};
\node at (-2.5,-1) {\large $\mathcal{E}\setminus \mathcal{D}$};
\node at (1.1,-2.5) {\large $\mathcal{D}$};

\end{scope}
\end{tikzpicture}
\caption{
 Two possible subdivisions of the set $\mathcal{E}$ when $\mathcal{V}$ is of type one. Left and right images correspond to  the cases  $\pro{C}{V_2-V_1}<\det V_1$ and $\pro{C}{V_2-V_1}>\det V_1$, respectively. }
\label{fig:div2}
\end{figure}
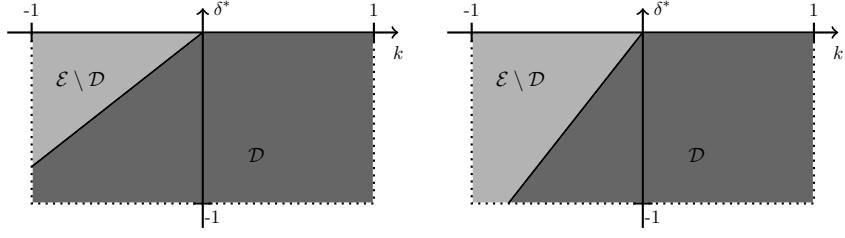In fact, $f_1$ and $f_2$ are real-valued functions, such that
\[
f_1\geq 0, \quad \mbox{and} \quad \kers f_1 = \{ V \in \msym \,| \, \pro{C}{V-V_2}\leq 0, \leq \det V, \,  0\leq \hbar(V)\},
\]
meanwhile,
\[
f_2\geq 0, \quad \mbox{and} \quad\kers f_2 \subset\{ V \in \msym \,| \,  0\leq \pro{C}{V-V_1},\, 0\leq \hbar(V)\}.
\] 

The latter equations imply that $f^{pp}$ is such that $f^{pp}\geq 0$ and 
\begin{equation}\label{kertype1}
 \kers f^{pp} \subset\{ V \in \msym \,| \,  \pro{C}{V-V_2}\leq 0\leq \pro{C}{V-V_1}, \ 0\leq \det V,\ 0\leq \hbar(V)\}.
\end{equation}
That is \cref{lem9:itemb2}, and the proof is complete.

 \end{proof}

\section{Proofs of main results}\label{sec:8}

This section concerns the proofs of Theorems \ref{T3} and \ref{T2}. 

\subsection{Proof of \cref{T3}. }

\begin{proof}
We consider three cases:
\begin{enumerate}[label=(\alph*)]
\item If all the wells in ${\mathcal U}$ are pairwise compatible Bhattacharya~\cite{B} p. 231, proved that 
$Q^e(\mathcal{U}) =L^e(\mathcal{U})=C(\mathcal{U})$,  and there is  nothing to prove. 
\item We assume that ${\mathcal U}$ is type two. Then, there exist a rank-one compatible pair of wells.  Without loss of 
generality, we assume that 
\begin{equation*}
\det(U_1-U_2)=0, \det(U_1-U_3)\leq 0, \text{ and } \det(U_2-U_3)>0.
\end{equation*}
Hence, by \cref{origin}, there exists $U_0\in C(\mathcal{U})$ such that $\det (U_2-U_0)=\det (U_3-U_0)=0$. Moreover, $U_0\in C(\{U_1,U_2\})$ and  $\det(U_1-U_0)=0$. In this case, the translated set $\mathcal{V}=\mathcal{U}-U_0$ satisfies that $\Aff(\mathcal{V})\subset \msym$ is a two-dimensional subspace, where $\det(V_i)=0$ for $i=1,2,3$. Hence, by the first part of \cref{eqboundary} (namely \cref{kerCP}), we conclude that if $f$ is given as in \cref{eq:fc}, then $f^{pp}$ is a non negative polyconvex function such that $\mathcal{V}\subset \kers f^{pp}$. 
Therefore, by \eqref{hullscontentions},
\[
L^e(\mathcal{V})\subset Q^e(\mathcal{V})\subset \kers f^{pp}\cap C(\mathcal{V})=\left\{V\in \msym\,\middle | \, V\in C(\mathcal{V}),\, 0 \leq \det V \right\}=C(\mathcal{V})\cap \overline{\mathcal{C}_0},
\] 
where the notation $\overline{\mathcal{C}_0}$ stands for the closure of the incompatible cone at $0$. Now, by \cref{P:lamconv}, we have that 
$L^e({\mathcal V})=C(\{0,V_2\})\cup C(\{0,V_1,V_3\})= C({\mathcal V})\cap \overline{\mathcal{C}_0}$, and the proof follows.

\item We assume that there is only one compatible (in fact rank-one compatible) pair in ${\mathcal U}$. In this case, ${\mathcal U}$ is of type one. Without loss of generality, we have that 
\begin{equation}\label{eq:T3-3}
\det(U_2-U_1)>0, \det(U_3-U_1)> 0 \text{ and } \det(U_3-U_2)=0.
\end{equation}
Now we claim that there exists $U_0\notin C({\mathcal U})$ such that 
\begin{equation}\label{eq:T3-4}
\det(U_0-U_1)=\det(U_0-U_2)=\det(U_0-U_3)=0.
\end{equation}


This claim follows by noticing that
$$
t\mapsto\det(tU_3+(1-t)U_2-U_1)=\det(U_2-U_1)+ t\pro{\mathcal{S}(U_3-U_2)}{U_2-U_1}
$$
 is linear in $t$ and $\pro{\mathcal{S}(U_3-U_2)}{U_2-U_1}\neq 0$, 
since $\Aff(\mathcal{U})$ is of codimension one. Thus, there exists $t_0\in\re$ such that \eqref{eq:T3-4} holds. Moreover, $t_0\notin [0,1]$, and $U_0=t_0U_3+(1-t_0)U_2\notin C({\mathcal U})$ by \eqref{eq:T3-3}, see \cref{fig:1}{\color{DarkOrchid}.(d)}.
Proceeding as before, we define $\mathcal{V}=\{V_1,\, V_2,\, V_3\}=\mathcal{U}-U_0$ where $V_i=U_i-U_0$  and $\det V_i=0$ for $i\in\{1,\,2,\,3\}$. Hence, by \cref{kerCP} and the first part of \cref{eqboundary}, we conclude that if $f$ is given as in \cref{eq:fc}, then $f^{pp}$ is a non negative polyconvex function such that $\mathcal{V}\subset \kers f^{pp}$, and
\begin{equation*}
L^e(\mathcal{V})\subset Q^e(\mathcal{V})\subset\kers f^{pp}\cap C({\mathcal V}) \subset  C({\mathcal V})\cap \overline{{\mathcal C}_{0}},
\end{equation*}
but 
$
C({\mathcal V})\cap \overline{{\mathcal C}_{0}}= \{V_3\}\cup C(\{V_1,V_2\})=L^e(\mathcal{V})
$ by \cref{P:lamconv}{\color{DarkOrchid}.(a)}, and the result follows.  
\end{enumerate}
\end{proof}

\begin{remark}\label{rmk5}
In the proof of \cref{T2} we are interested in $\Ker f^{pp}\cap C(\mathcal{V})$. Under the assumptions of \cref{eqboundary} item 2, we claim that the conditions $V\in  C(\mathcal{V})$ and $0\le \hbar(V) $ are equivalent to $\theta_1\det V_1 \le \det V$ for some $\theta_1\in [0,1]$. Assume the claim for a moment. We notice that the latter equation imply that if $\mathcal{V}$ is of type one or type two, then $0 < \det V_1$ and $0\le \det V$, or  $\det V_1<0$ and $\det V_1 < \det V$, respectively. Hence, from \eqref{kertype1} and \eqref{kertype2}, we get  
\[
\Ker f^{pp}\cap C(\mathcal{V})\subset \{ M \in C(\mathcal{V})| \pro{C}{e(M)-V_2} \leq 0  \leq \pro{C}{e(M)-V_1}, 0\leq \hbar(M)\},
\] 
for both cases. 
Now, we prove the claim. By assumption $\pro{C}{V_1-V_2}\neq 0$, thus from $0\le \hbar(V)$ we get that 
\begin{equation}\label{alth}
\frac{\pro{C}{V-V_2}}{\pro{C}{V_1-V_2}}\det V_1 \leq \det V, 
\end{equation}
Now,  since $V\in C(\mathcal{V})$, then $V=\theta_1 V_1+\theta_2 V_2 +\theta_3 V_3$ for some $\theta_1,\theta_2,\theta_3\in [0,1]$ such that 
$\theta_1+\theta_2+\theta_3=1$. Thus, $\pro{C}{V-V_2}=\pro{C}{\theta_1V_1-(1-\theta_2-\theta_3)V_2}=\theta_1\pro{C}{V_1-V_2}$ due to  $\pro{C}{V_2}=\pro{C}{V_3}$. Hence, 
\cref{alth} implies one side of the claim, and $\theta_1$ is the volume fraction of $V_1$. The reverse implication follows simple algebraic manipulations. 
\end{remark}


\subsection{Proof of \cref{T2}. }

\begin{proof}
We notice that under the assumptions of \cref{T2}, \cref{origin} guarantees the existence of $U_{0}\in C(\mathcal{U})$ such that $\det(U_2-U_{0})=\det(U_3-U_{0})=0$ and $ \det(U_1-U_{0})\neq 0$. Moreover, the translated set $\mathcal{V}=\{V_1,V_2,V_3\}:=\mathcal{U}-U_{0}$ satisfies that $0\in C(\mathcal{V})$, $\det(V_2)=\det(V_3)=0$, and $\Aff(\mathcal{V})$ is a two dimensional subspace.  Hence, by \cref{eqboundary}, we obtain explicit bound on $Q^e(\mathcal{U})$ by intersecting  $C(\mathcal{U})$ with $\kers f^{pp}$ for admissible matrices $C$.  Then by \cref{trans} and \cref{hullscontentions}, 
\begin{equation} \label{eq:qqq}
Q^e(\mathcal{U})\subset \left(\kers f^{pp}+U_{0}\right)\cap C(\mathcal{U}),
\end{equation}
where $\kers f^{pp}+U_{0}$ stands for the translated set $\{U_0+V\,|\, V\in \kers f^{pp}\}$.
By \cref{rem:det<}, if $Q\in\msym$ is orthogonal to $\Aff({\mathcal{V}})$, then $\det Q<0$, and due to \cref{cor:origin},
\begin{equation}\label{eq:vs}
V_1 = \pv{\xi}{\eta}, \quad V_2 = \pva{\gamma}, \quad \mbox{and} \quad V_3 = \pvn{\zeta},
\end{equation}
for some linearly independent $a,n\in S^1$ and constants $\xi,\eta,\gamma,\zeta\in \re$ that satisfy \eqref{coefU}. 

Next, in order to apply \cref{eqboundary}, we construct a proper matrix $C$. To this end,  we denote by $\tilde{V}$ the $\mathbb{R}^3$ representation of $V\in\msym$ under the isomorphism  \eqref{isomorphism}, and let
\[
\tilde{C}_0=\tilde{Q}\times (\tilde{V}_2-\tilde{V}_3)
 \quad \mbox{where} \quad 
 \tilde{Q}= \frac{(\tilde{V}_2-\tilde{V}_1)\times (\tilde{V}_3-\tilde{V}_1)}{\|(\tilde{V}_2-\tilde{V}_1)\times (\tilde{V}_3-\tilde{V}_1)\|}.  
\]
By construction, the matrix $C_0$ belongs to the affine space $\Aff(\mathcal{V})$, and it is perpendicular to $V_3-V_2$.
After some direct algebraic manipulations, we get 
\begin{equation}\label{C_0}
\tilde{C}_0=
\frac{(\tilde{V}_2-\tilde{V}_1)\pro{\tilde{V}_3-\tilde{V}_1}{\tilde{V}_3-\tilde{V}_2}+(\tilde{V}_3-\tilde{V}_1)\pro{\tilde{V}_2-\tilde{V}_1}{\tilde{V}_2-\tilde{V}_3}}
{|P|},
\end{equation}
where
\begin{equation*}
P^2:=\|(\tilde{V}_2-\tilde{V}_1)\times (\tilde{V}_3-\tilde{V}_1)\|^2=\|\tilde{V}_2-\tilde{V}_1\|^2\|\tilde{V}_3-\tilde{V}_1)\|^2-\pro{\tilde{V}_2-\tilde{V}_1}{\tilde{V}_3-\tilde{V}_1}^2.
\end{equation*}
 Notice that, since the isomorphism \eqref{isomorphism}  preserves the inner product, we can drop the tildes out in \cref{C_0} to get the corresponding expression in $\msym$. 
In terms of the coordinates given by \eqref{eq:vs}, we get that
\[
P^2= (\eta\gamma+\xi\zeta-\gamma\zeta)^2(1-(a\cdot n)^4), 
\]
and
\begin{equation*} 
C_0=-\left(\frac{\left[\pv{\left(\zeta-\gamma(a\cdot n)^2\right)}{\left(\gamma-\zeta(a\cdot n)^2\right)}\right]}{\sqrt{1-(a\cdot n)^4}}\right) \sgn(\eta\gamma+\xi\zeta-\gamma\zeta).
\end{equation*}
By \cref{coefU} and \cref{cor:origin}, 
we get 
\[
\sgn(\eta\gamma+\xi\zeta-\gamma\zeta)=\begin{cases}-1 &\mathcal{U} \mbox{ is type one,}\\1 &\mathcal{U} \mbox{ is type two.}\\
\end{cases}
\] 
Thus, 
\[
\det C_0 = 
\frac{\gamma \zeta(1+(a\cdot n)^2)^2-(\gamma+\zeta)^2(a\cdot n)^2}{1+(a\cdot n)^2}.
\] 
First, we  assume $\mathcal{U}$ is type two. Then, $\mathcal{V}$ is also type two, and we  let  
\[
C=\frac{1}{\sqrt{1-(a\cdot n)^4}}C_0.
\] 
In this case, $\det C<0$, and $\Pi_C\cap \Aff(\mathcal{V})$ is a line, parallel to $V_2-V_3$. Moreover, by an explicit computation and \eqref{coefU}, we get
\begin{equation*}
  \pro{C}{V_{2}}=\pro{C}{V_{3}}= -\zeta \gamma>0\geq -(\xi\zeta+\eta\gamma)=\pro{C}{V_1}.
\end{equation*}
Therefore, the set $\mathcal{V}$ and the matrix $C$ meet necessary conditions to apply \cref{eqboundary} part \ref{lem9:itemb1} and \cref{rmk5}. Thus,  we have
\begin{equation*}
\kers f^{pp}=\{ V \in \msym \,| \, \pro{C}{V-V_2} \leq 0  \leq \pro{C}{V-V_1}, \ \  0\leq \hbar(V)\}.
\end{equation*}
In view of \cref{eq:qqq}, we are interested in the set $\kers f^{pp}\cap C(\mathcal{V})$. We notice that, the sets $\{ V \in \Aff(\mathcal{V})\,| \, \pro{C}{V-V_2}=0\}$ and $\{ V \in \Aff(\mathcal{V})\,| \, \pro{C}{V-V_3}=0\}$  are supporting lines of $C(\mathcal{V})$ on $V_2$ and $V_1$ respectively. Hence, 
\begin{equation}\label{eq:1erbound}
Q^e(\mathcal{V})\subset\kers f^{pp}\cap C(\mathcal{V}) = \{ V \in C(\mathcal{V}) \,| \ \  0\leq \hbar(V)\},
\end{equation}
with
\begin{equation}\label{eq:hbar}
    \hbar(V) = \pro{C}{V-V_2}\det V_1-\pro{C}{V_1-V_2}\det V.
\end{equation}
Now, \cref{T2} for ${\mathcal U}$ of type two follows from the assignment $V_i=U_i +U_0$ and $V=U+U_0$. 
\medskip 

Second, we  assume $\mathcal{U}$ is type one. 
Then, $\mathcal{V}$ is also type one.  If we let $Q\in \msym$ be the  orthogonal direction to $\Aff(\mathcal{V})$, then we define $C(t)\in\msym$ for each $t\in(-\pi,\pi]$ as  
\begin{equation*}
C(t)=Q\cos t +C\sin t,
\end{equation*}
where $C$ is chosen as before  in terms of $C_0$.   In this case, we can readily see that 
\[
\pro{C(t)}{V_{2}}=\pro{C(t)}{V_{3}}= \zeta \gamma\sin t,\quad \mbox{and}\quad
\pro{C(t)}{V_1}=(\xi\zeta+\eta\gamma)\sin t .
\]
Thus, due to the continuity of $t\mapsto \det C(t)$ and \cref{rem:det<}, there exists $t_0>0$ small enough such that $\det C(t_0)<0$ and $\pro{C(t_0)}{V_1}\leq 0< \pro{C(t_0)}{V_3}=\pro{C(t_0)}{V_2}$. Hence, by  \cref{eqboundary}~\ref{lem9:itemb2} and \cref{rmk5}, we have that 
\[\kers f^{pp}=\{ V \in \msym \,| \, \pro{C(t_0)}{V-V_2} \leq 0  \leq \pro{C(t_0)}{V-V_1}, \ \  0\leq \kappa(V)\},
\]
where $\kappa(V) = \pro{C(t_0)}{V-V_2}\det V_1-\pro{C(t_0)}{V_1-V_2}\det V.$ By a similar argument as before, we conclude that
\begin{equation*}
Q^e(\mathcal{V})\subset \kers f^{pp}\cap C(\mathcal{V})=\{ V \in C(\mathcal{V}) \,| \,   0\leq \kappa(V)\}.
\end{equation*}
Next, since $V\in C(\mathcal{V})$, it follows that $\kappa(V)=\sin t_0\, \hbar(V)$ where $\hbar$ is defined by \cref{eq:hbar}. Thus, we conclude that 
\begin{equation}\label{eq:2ndbound}
Q^e(\mathcal{V})\subset\{ V \in C(\mathcal{V}) \,| \  0\leq \hbar(V)\},
\end{equation}
and \cref{T2} follows from \cref{eq:qqq}. Therefore, we finish the proof of our theorem. 
\end{proof}
\begin{remark}
The component-wise representation of \cref{eq:1erbound} and \cref{eq:2ndbound} is straight forward.  Indeed, if $V\in \Aff(\mathcal{V})$, then $V=\pva{x}+ \pvn{y}$ for some  real $x,y$, and $\hbar(V)=|n\times a|^2h(x,y)$, where
\begin{equation*}
h(x,y)=x y (\eta \gamma + (\xi - \gamma) \zeta) - \xi \eta (\gamma y + (x-\gamma) \zeta ).
\end{equation*}
Therefore, \cref{eq:1erbound} is rewritten as
\begin{equation*}
Q^e(\mathcal{V}) \subset \left\{V \in C(\mathcal{V}) \, |\,V=\pv{x}{y},\, \  0\leq h(x,y)  \right\}.
\end{equation*}
Meanwhile, \cref{eq:2ndbound} is rewritten as
\begin{equation*}
Q^e(\mathcal{V}) \subset \left\{V \in C(\mathcal{V}) \, |\,U=\pv{x}{y},\, \  h(x,y)\leq 0  \right\}.
\end{equation*}
\end{remark}

\section[Optimal bounds on \texorpdfstring{$Q^e(\mathcal{U})$}{} by  quadratic polyconvex functions]{Optimality of the bound for \texorpdfstring{$Q^e(\mathcal{U})$}{} by symmetric quadratic polyconvex functions}\label{sec:9}

As a final remark, we comment briefly about the optimality of the bounds in \cref{T2}. By the characterization of symmetric quadratic polyconvex functions given by Boussaid {\em et. al} \cite{Anja19} (proposition 4.5 in page 435), we can prove that our outer bound for $Q^e(\mathcal{U})$ is optimal when we restrict the analysis to quadratic functions.  
\medskip 

Let $\nu$ be a homogeneous Young measure limit of linear strains supported on $\mathcal{V}=\mathcal{U}-U_0$, {\em i.e.} $\nu= \theta_1\delta_{V_1}+\theta_2\delta_{V_2}+\theta_3\delta_{V_3}$ 
with barycenter  $V=\theta_1V_1+\theta_2V_2+\theta_3V_3$. Since quasiconvex functions preserves Jensen's inequality almost everywhere in $\Omega$, see for example lemma 1.1 in \cite{BFJK} p. 850, we have that  
\[
f\left(\int_{\msym}Ad\,\nu(A)\right)\leq \int_{\msym}f(A)d\,\nu(A).
\]
 Now, by proposition 4.5 in \cite{Anja19}, $f:\msym \rightarrow \re$ is symmetric quadratic polyconvex if and only if it has the form  $f(\cdot )=g(\cdot ) - \alpha \det(\cdot )$  for some convex function $g:\msym \rightarrow \re$ and $\alpha>0$. Therefore $f$ is also symmetric quasiconvex. Hence, 
 \[
g(V) -\alpha \det(V) \leq \sum_{i=1}^3 \theta_i g(V_i) -\alpha  \sum_{i=1}^3 \theta_i \det(V_i).
\]
Since, under the assumptions of \cref{T2} and the definition of $U_0$,  $\theta_2\det V_2=\theta_3\det V_3=0$. Therefore, 
\begin{equation}\label{YMB2}
\det(V) - \theta_1 \det(V_1) \geq \sup \left\{ 
\frac{1}{\alpha}\left(g(V)-\sum_{i=1}^3 \theta_i g(V_i)\right)\ \middle |\  g \ \mbox{is convex}, \alpha>0 \right \},
\end{equation}
for every $V\in Q^e(\mathcal{V})$. The supremum in the right hand side of \cref{YMB2} is attained by the convex function $g=0$. Hence, $\det(V) - \theta_1 \det(V_1)\geq 0$ and by \cref{rmk5}, we recover $V\in C({\mathcal V})$ and $0\le \hbar(V)$. Therefore, we obtain an outer bound for  $Q^e(\mathcal{V})$ given by
\begin{equation}\label{hiddenh}
\left\{V=\theta_1V_1+\theta_2V_2+\theta_3V_3\ \middle | \ \det(V) \geq \theta_1 \det(V_1),\, (\theta_1,\theta_2,\theta_3)\in [0,1]^3,\,  \sum_{i=1}^3 \theta_i=1\right\} 
\end{equation}
The  bound in \cref{hiddenh} is equivalent to the outer in \cref{T2}.


\section{Acknowledgements}
The authors kindly acknowledges the support given by UNAM PAPPIT--IN106118 grant. 
LM acknowledge the support of CONACyT grant--590176 during his graduate studies 
and AC was also partially founded by CONACYT CB-2016-01-284451 grant.

 \section*{Conflict of interest}
 The authors declare that they have no conflict of interest.

\bibliographystyle{abbrv}      
\bibliography{MyRef2}   

\end{document}